\documentclass[envcountsame,envcountresetchap]{svmult}
\usepackage[bottom]{footmisc}
\def\arrows_for_Sasha#1#2{ 
\setlength{\unitlength}{.8cm}
\begin{picture}(1.5,0.00)(-0.25,0.00)
\put(0,.25){\vector(1,0){1}}
\put(1,0){\vector(-1,0){1}}
\put(.5,-.15){\makebox(0,0)[t]{\scriptsize$#1$}}
\put(.5,.4){\makebox(0,0)[b]{\scriptsize$#2$}}
\end{picture}
}
\def\itIB{{\it IB}} \def\yyy{}
\def\cdotvetsi{\hskip.1em{\mbox {\tiny \raisebox{.2em}{$\bullet$}}}\hskip.1em}
\def\itLieB{{\it LieB}}
\def\calP{{\cal P}}

\def\calQ{{\cal Q}}

\def\sfP{{\sf P}}
\def\id{{\mbox{1 \hskip -8pt 1}}}
\def\PROP{{\small\sc PROP}}
\def\hPROP{$\frac12${\small\sc PROP}}
\def\catPROP{{\tt PROP}}
\def\cathPROP{{\textstyle\frac12}{\tt PROP}}
\def\cathPROPsoft{{\frac12}{\tt PROP}}
\def\catOper{{\tt Oper}} \def\catProper{\tt Proper}
\def\catDiop{{\tt diOp}}
\def\ot{\otimes}
\def\End{{\it {\cal E}\hskip -1pt  nd}}
\def\Ass{{\it {\cal A}  ss}}
\def\Lie{{\it {\cal L}  ie}}
\def\Com{{\it {\cal C}  om}}
\def\jcirc#1{{\hskip 1mm {}_{#1}\circ \hskip .2mm}}
\def\jcircsoft#1#2{{\hskip 1mm {}_{#1}\circ_{#2}}}
\def\Hom{{\it Hom}}
\def\otexp#1#2{{#1}^{\otimes #2}}
\def\papert{\partial_{\it pert}}
\def\sfB{{\sf B}}
\def\sfIB{{\sf IB}}
\def\sfLieB{{\sf LieB}}
\def\br{{\rm br}}
\def\span{{\it Span}}
\def\catbCol{{\tt bCol}}
\def\sfhB{{\mbox{$\frac12$\sf B}}}
\def\sfhb{{\mbox{$\frac12$\sf b}}}
\def\sfhlieb{{\mbox{$\frac12$\sf lieb}}}

\def\ithLieB{{\mbox{$\frac12$\it LieB}}}
\def\sfhLieB{{\mbox{$\frac12$\sf LieB}}}

\def\sfs{{\sf s}}
\def\sfP{{\sf P}}
\def\pa{\partial}
\def\sft{{\sf t}}
\def\gns{{\rm gen}} \def\comp{{\rm cmp}} \def\grad{{\rm grad}}
\def\pth{{\rm pth}} \def\gencomp{{\rm c+g}}
\def\Associative{\Ass}
\def\BOX{\raisebox{1.2mm}%
         {\hskip .5mm $\fbox{\hphantom{\hglue .01mm}}$\hskip .5mm}}
\def\Box{\BOX}

\def\dvojiteypsilon{{
\unitlength=.3pt
\begin{picture}(24.00,30.00)(0.00,3.00)
\put(10.00,20.00){\line(0,-1){10.00}}
\bezier{20}(10.00,10.00)(15,5)(20.00,0.00)
\bezier{20}(10.00,10.00)(5,5)(0.00,0.00)
\bezier{20}(10.00,20.00)(15,25)(20.00,30.00)
\bezier{20}(0.00,30.00)(5,25)(10.00,20.00)
\end{picture}}}

\def\dvojiteypsilonvetsi#1#2#3#4{{
\unitlength=.46pt
\begin{picture}(24.00,30.00)(0.00,3.00)
\put(10.00,20.00){\line(0,-1){10.00}}
\bezier{20}(10.00,10.00)(15,5)(20.00,0.00)
\bezier{20}(10.00,10.00)(5,5)(0.00,0.00)
\bezier{20}(10.00,20.00)(15,25)(20.00,30.00)
\bezier{20}(0.00,30.00)(5,25)(10.00,20.00)
\put(0,-5){\makebox(0,0)[t]{\scriptsize $#1$}}
\put(20,-5){\makebox(0,0)[t]{\scriptsize $#2$}}
\put(0,35){\makebox(0,0)[b]{\scriptsize $#3$}}
\put(20,35){\makebox(0,0)[b]{\scriptsize $#4$}}
\end{picture}}}

\def\pravaplastev{
\unitlength=.35pt
\begin{picture}(38,30)(-4,4)
\put(0.00,0.00){\line(0,1){10}}
\put(20.00,0.00){\line(0,1){10}}
\put(10.00,20.00){\line(0,1){10}}
\put(30.00,20.00){\line(0,1){10}}
\put(0,10){\bezier{20}(0.00,0.00)(5.00,5.00)(10.00,10.00)}
\put(20,10){\bezier{20}(0.00,0.00)(5.00,5.00)(10.00,10.00)}
\put(20,10){\bezier{20}(0.00,0.00)(-5.00,5.00)(-10.00,10.00)}
\end{picture}
}

\def\pravaplastevvetsi#1#2#3#4{
\unitlength=.45pt
\begin{picture}(38,30)(-4,4)
\put(0.00,0.00){\line(0,1){10}}
\put(20.00,0.00){\line(0,1){10}}
\put(10.00,20.00){\line(0,1){10}}
\put(30.00,20.00){\line(0,1){10}}
\put(0,10){\bezier{20}(0.00,0.00)(5.00,5.00)(10.00,10.00)}
\put(20,10){\bezier{20}(0.00,0.00)(5.00,5.00)(10.00,10.00)}
\put(20,10){\bezier{20}(0.00,0.00)(-5.00,5.00)(-10.00,10.00)}
\put(0,-5){\makebox(0,0)[t]{\scriptsize $#1$}}
\put(20,-5){\makebox(0,0)[t]{\scriptsize $#2$}}
\put(10,35){\makebox(0,0)[b]{\scriptsize $#3$}}
\put(30,35){\makebox(0,0)[b]{\scriptsize $#4$}}
\end{picture}
}

\def\levaplastev{
\unitlength=.35pt
\begin{picture}(38,30)(-34,4)
\put(0.00,0.00){\line(0,1){10}}
\put(-20.00,0.00){\line(0,1){10}}
\put(-10.00,20.00){\line(0,1){10}}
\put(-30.00,20.00){\line(0,1){10}}
\put(0,10){\bezier{20}(0.00,0.00)(-5.00,5.00)(-10.00,10.00)}
\put(-20,10){\bezier{20}(0.00,0.00)(-5.00,5.00)(-10.00,10.00)}
\put(-20,10){\bezier{20}(0.00,0.00)(5.00,5.00)(10.00,10.00)}
\end{picture}
}

\def\levaplastevvetsi#1#2#3#4{
\unitlength=.45pt
\begin{picture}(38,30)(-34,4)
\put(0.00,0.00){\line(0,1){10}}
\put(-20.00,0.00){\line(0,1){10}}
\put(-10.00,20.00){\line(0,1){10}}
\put(-30.00,20.00){\line(0,1){10}}
\put(0,10){\bezier{20}(0.00,0.00)(-5.00,5.00)(-10.00,10.00)}
\put(-20,10){\bezier{20}(0.00,0.00)(-5.00,5.00)(-10.00,10.00)}
\put(-20,10){\bezier{20}(0.00,0.00)(5.00,5.00)(10.00,10.00)}
\put(0,-5){\makebox(0,0)[t]{\scriptsize $#2$}}
\put(-20,-5){\makebox(0,0)[t]{\scriptsize $#1$}}
\put(-10,35){\makebox(0,0)[b]{\scriptsize $#4$}}
\put(-30,35){\makebox(0,0)[b]{\scriptsize $#3$}}
\end{picture}
}

\def\jitka{
\unitlength=.3pt
\begin{picture}(40,20)(-2,0)
\put(0,0){\line(0,1){10}}
\put(10,20){\line(0,1){10}}
\put(30,20){\line(0,1){10}}
\put(10,0){\bezier{40}(0.00,0.00)(10.00,10.00)(20.00,20.00)}
\put(10,0){\bezier{40}(0.00,20.00)(10.00,10.00)(20.00,0.00)}
\put(0,10){\bezier{20}(0.00,0.00)(5.00,5.00)(10.00,10.00)}
\end{picture}
}

\def\jitkainv{
\unitlength=.3pt
\begin{picture}(40,20)(-34,0)
\put(0,0){\line(0,1){10}}
\put(-10,20){\line(0,1){10}}
\put(-30,20){\line(0,1){10}}
\put(-10,0){\bezier{40}(0.00,0.00)(-10.00,10.00)(-20.00,20.00)}
\put(-10,0){\bezier{40}(0.00,20.00)(-10.00,10.00)(-20.00,0.00)}
\put(0,10){\bezier{20}(0.00,0.00)(-5.00,5.00)(-10.00,10.00)}
\end{picture}
}

\def\anna{
\unitlength=.32pt
\begin{picture}(40,20)(-2,0)
\put(0,0){\line(0,1){10}}
\put(10,0){\line(0,1){30}}
\put(20,0){\line(0,1){10}}
\put(30,20){\line(0,1){10}}
\put(0,10){\bezier{20}(0.00,0.00)(5.00,5.00)(10.00,10.00)}
\put(20,10){\bezier{20}(0.00,0.00)(5.00,5.00)(10.00,10.00)}
\put(20,10){\bezier{20}(0.00,0.00)(-5.00,5.00)(-10.00,10.00)}
\end{picture}
}

\def\annainv{
\unitlength=.32pt
\begin{picture}(40,20)(-34,0)
\put(0,0){\line(0,1){10}}
\put(-10,0){\line(0,1){30}}
\put(-20,0){\line(0,1){10}}
\put(-30,20){\line(0,1){10}}
\put(0,10){\bezier{20}(0.00,0.00)(-5.00,5.00)(-10.00,10.00)}
\put(-20,10){\bezier{20}(0.00,0.00)(-5.00,5.00)(-10.00,10.00)}
\put(-20,10){\bezier{20}(0.00,0.00)(5.00,5.00)(10.00,10.00)}
\end{picture}
}

\def\motylek{{
\unitlength=.3pt
\begin{picture}(66.00,60.00)(-3.00,20.00)
\put(10.00,50.00){\line(0,1){0.00}}
\put(50.00,10.00){\line(0,-1){10.00}}
\put(10.00,10.00){\line(0,-1){10.00}}
\put(50.00,60.00){\line(0,-1){10.00}}
\put(10.00,60.00){\line(0,-1){10.00}}
\put(60.00,40.00){\line(0,-1){20.00}}
\put(0.00,40.00){\line(0,-1){20.00}}
\bezier{20}(50.00,10.00)(55.00,15.00)(60.00,20.00)
\bezier{20}(0.00,20.00)(5.00,15.00)(10.00,10.00)
\bezier{20}(50.00,50.00)(55.00,45.00)(60.00,40.00)
\bezier{20}(10.00,10.00)(15,15)(22,22)
\put(28,28){\bezier{20}(10.00,10.00)(15,15)(22,22)}
\put(10.00,50.00){\line(1,-1){40.00}}
\bezier{20}(10.00,50.00)(5.00,45.00)(0.00,40.00)
\end{picture}}}

\def\jednadva{{
\unitlength=.4pt
\begin{picture}(24.00,20.00)(-2.00,0.00)
\bezier{20}(10.00,10.00)(15.00,5.00)(20.00,0.00)
\bezier{20}(10.00,10.00)(5.00,5.00)(0.00,0.00)
\put(10.00,20.00){\line(0,-1){10.00}}
\end{picture}}
}

\def\jednactyri{{
\unitlength=.05pt
\begin{picture}(176.00,160.00)(-8.00,0.00)
\put(80.00,100.00){\line(0,1){60.00}}
\bezier{20}(80.00,80.00)(100.00,30.00)(110.00,0.00)
\bezier{20}(80.00,80.00)(60.00,30.00)(50.00,0.00)
\bezier{20}(80.00,80.00)(120.00,40.00)(160.00,0.00)
\bezier{20}(80.00,80.00)(40.00,40.00)(0.00,0.00)
\put(80.00,100.00){\line(0,-1){20.00}}
\end{picture}}
}

\def\ctyrijedna{{
\unitlength=.05pt
\begin{picture}(176.00,160.00)(-8.00,-160.00)
\put(80.00,-100.00){\line(0,-1){60.00}}
\bezier{20}(80.00,-80.00)(100.00,-30.00)(110.00,0.00)
\bezier{20}(80.00,-80.00)(60.00,-30.00)(50.00,0.00)
\bezier{20}(80.00,-80.00)(120.00,-40.00)(160.00,0.00)
\bezier{20}(80.00,-80.00)(40.00,-40.00)(0.00,0.00)
\put(80.00,-80.00){\line(0,-1){40.00}}
\end{picture}}
}

\def\dvajedna{{
\unitlength=.4pt
\begin{picture}(24.00,20.00)(-2.00,0.00)
\put(10.00,10.00){\line(0,-1){10.00}}
\bezier{20}(10.00,10.00)(15.00,15.00)(20.00,20.00)
\bezier{20}(0.00,20.00)(5.00,15.00)(10.00,10.00)
\end{picture}}
}

\def\dvadva{{
\unitlength=.8pt
\begin{picture}(12.00,10.00)(-1.00,0.00)
\bezier{30}(0.00,0.00)(5.00,5.00)(10.00,10.00)
\bezier{30}(0.00,10.00)(5.00,5.00)(10.00,0.00)
\end{picture}}
}

\def\jednatri{{
\unitlength=.4pt
\begin{picture}(24.00,20.00)(-2.00,0.00)
\bezier{20}(10.00,10.00)(15.00,5.00)(20.00,0.00)
\bezier{20}(10.00,10.00)(5.00,5.00)(0.00,0.00)
\put(10.00,20.00){\line(0,-1){20.00}}
\end{picture}}
}

\def\trijedna{{
\unitlength=.4pt
\begin{picture}(24.00,20.00)(-2.00,-20.00)
\bezier{20}(10.00,-10.00)(15.00,-5.00)(20.00,0.00)
\bezier{20}(10.00,-10.00)(5.00,-5.00)(0.00,0.00)
\put(10.00,-20.00){\line(0,1){20.00}}
\end{picture}}
}

\def\dvatri{{
\unitlength=0.4pt
\begin{picture}(24.00,20.00)(-2.00,0.00)
\put(10.00,10.00){\line(0,-1){10.00}}
\bezier{30}(0.00,0.00)(10.00,10.00)(20.00,20.00)
\bezier{30}(0.00,20.00)(10.00,10.00)(20.00,0.00)
\end{picture}}
}

\def\tridva{{
\unitlength=.4pt
\begin{picture}(24.00,20.00)(-2.00,-20.00)
\put(10.00,-10.00){\line(0,1){10.00}}
\bezier{30}(0.00,0.00)(10.00,-10.00)(20.00,-20.00)
\bezier{30}(0.00,-20.00)(10.00,-10.00)(20.00,0.00)
\end{picture}}
}

\def\tritri{{
\unitlength=.4pt
\begin{picture}(24.00,20.00)(-2.00,0.00)
\bezier{30}(0.00,0.00)(10.00,10.00)(20.00,20.00)
\bezier{30}(0.00,20.00)(10.00,10.00)(20.00,0.00)
\put(10.00,20.00){\line(0,-1){20.00}}
\end{picture}}
}

\def\dvactyri{{
\unitlength=.1pt
\begin{picture}(96.00,80.00)(-8.00,0.00)
\bezier{30}(0.00,0.00)(40.00,40.00)(80.00,80.00)
\bezier{30}(0.00,80.00)(40.00,40.00)(80.00,0.00)
\bezier{20}(40.00,40.00)(48.50,14.00)(53.50,0.00)
\bezier{20}(40.00,40.00)(32.50,16.00)(26.00,0.00)
\end{picture}}
}

\def\dvaZbbbZb{{
\unitlength=0.1pt
\begin{picture}(96.00,80.00)(-8.00,0.00)
\bezier{10}(20.00,20.00)(30.00,10.00)(40.00,0.00)
\put(20.00,20.00){\line(0,-1){20.00}}
\bezier{30}(0.00,0.00)(40.00,40.00)(80.00,80.00)
\bezier{30}(0.00,80.00)(40.00,40.00)(80.00,0.00)
\end{picture}}
}

\def\dvabZbbbZ{{
\unitlength= 0.1pt
\begin{picture}(96.00,80.00)(-8.00,0.00)
\bezier{10}(60.00,20.00)(50.00,10.00)(40.00,0.00)
\put(60.00,20.00){\line(0,-1){20.00}}
\bezier{30}(0.00,0.00)(40.00,40.00)(80.00,80.00)
\bezier{30}(0.00,80.00)(40.00,40.00)(80.00,0.00)
\end{picture}}
}

\def\dvaZbbZbb{{
\unitlength=0.1pt
\begin{picture}(96.00,80.00)(-8.00,0.00)
\bezier{8}(11.00,9.00)(15.50,4.00)(21.00,0.00)
\put(40.00,40.00){\line(0,-1){40.00}}
\bezier{30}(0.00,0.00)(40.00,40.00)(80.00,80.00)
\bezier{30}(0.00,80.00)(40.00,40.00)(80.00,0.00)
\end{picture}}
}

\def\dvabbZbbZ{{
\unitlength=0.1pt
\begin{picture}(96.00,80.00)(-8.00,0.00)
\bezier{5}(70.00,10.00)(65.00,5.50)(60.00,0.00)
\put(40.00,40.00){\line(0,-1){40.00}}
\bezier{30}(0.00,0.00)(40.00,40.00)(80.00,80.00)
\bezier{30}(0.00,80.00)(40.00,40.00)(80.00,0.00)
\end{picture}}
}

\def\dvabZbbZb{{
\unitlength=0.1pt
\begin{picture}(96.00,80.00)(-8.00,0.00)
\bezier{10}(40.00,20.00)(50.00,10.00)(60.00,0.00)
\bezier{10}(40.00,20.00)(30.00,10.00)(20.00,0.00)
\put(40.00,40.00){\line(0,-1){20.00}}
\bezier{30}(0.00,0.00)(40.00,40.00)(80.00,80.00)
\bezier{30}(0.00,80.00)(40.00,40.00)(80.00,0.00)
\end{picture}}
}

\def\dvacarkatri{{
\unitlength=.2pt
\begin{picture}(40.00,50.00)(0.00,0.00)
\put(20.00,30.00){\line(0,-1){10.00}}
\put(20.00,0.00){\line(0,1){20}}
\bezier{20}(20.00,20.00)(30.00,10.00)(40.00,0.00)
\bezier{20}(20.00,20.00)(10.00,10.00)(0.00,0.00)
\bezier{20}(20.00,30.00)(30.00,40.00)(40.00,50.00)
\bezier{20}(0.00,50.00)(10.00,40.00)(20.00,30.00)
\end{picture}}
}

\def\gen#1#2{
\if #11
    \if #22 \jednadva \else \fi
\else
\fi
\if #12
    \if #22 \dvadva \else \fi
\else
\fi
\if #12
    \if #21 \dvajedna \else \fi
\else
\fi
\if #13
    \if #22 \tridva \else \fi
\else
\fi
\if #13
    \if #21 \trijedna \else \fi
\else
\fi
\if #12
    \if #23 \dvatri \else \fi
\else
\fi
\if #11
    \if #23 \jednatri \else \fi
\else
\fi
\if #11
    \if #24 \jednactyri \else \fi
\fi
\if #14
    \if #21 \ctyrijedna \else \fi
\fi
}

\def\bZbbZ{
{
\unitlength=.27pt
\begin{picture}(48.00,30.00)(-4,0.00)
\bezier{34}(20.00,20.00)(30.00,10.00)(40.00,0.00)
\bezier{34}(20.00,20.00)(10.00,10.00)(0.00,0.00)
\bezier{20}(30.00,10.00)(25.00,5.00)(20.00,0.00)
\put(20.00,30.00){\line(0,-1){10.00}}
\end{picture}} 
}

\def\ZbbZb{{
\unitlength=.27pt
\begin{picture}(48.00,30.00)(-4,0.00)
\bezier{34}(20.00,20.00)(30.00,10.00)(40.00,0.00)
\bezier{34}(20.00,20.00)(10.00,10.00)(0.00,0.00)
\bezier{20}(10.00,10.00)(15.00,5.00)(20.00,0.00)
\put(20.00,30.00){\line(0,-1){10.00}}
\end{picture}}}

\def\Jac#1#2#3{{
\unitlength=.4pt
\begin{picture}(48.00,30.00)(-4,0.00)
\bezier{34}(20.00,20.00)(30.00,10.00)(40.00,0.00)
\bezier{34}(20.00,20.00)(10.00,10.00)(0.00,0.00)
\bezier{20}(10.00,10.00)(15.00,5.00)(20.00,0.00)
\put(20.00,30.00){\line(0,-1){10.00}}
\put(0,-5){\makebox(0,0)[t]{\scriptsize $#1$}}
\put(20,-5){\makebox(0,0)[t]{\scriptsize $#2$}}
\put(40,-5){\makebox(0,0)[t]{\scriptsize $#3$}}
\end{picture}}}

\def\ZbbZbb{{
\unitlength=.07pt
\begin{picture}(192.00,120.00)(-16.00,0.00)
\bezier{20}(20.00,20.00)(30.00,10.00)(40.00,0.00)
\bezier{30}(80.00,80.00)(120.00,40.00)(160.00,0.00)
\bezier{30}(80.00,80.00)(40.00,40.00)(0.00,0.00)
\put(80.00,120.00){\line(0,-1){120.00}}
\end{picture}}
}

\def\bbZbbZ{{
\unitlength=0.07pt
\begin{picture}(192.00,120.00)(-16.00,0.00)
\bezier{20}(140.00,20.00)(130.00,10.00)(120.00,0.00)
\bezier{30}(80.00,80.00)(120.00,40.00)(160.00,0.00)
\bezier{30}(80.00,80.00)(40.00,40.00)(0.00,0.00)
\put(80.00,120.00){\line(0,-1){120.00}}
\end{picture}}
}

\def\bZbbZb{{
\unitlength=.07pt
\begin{picture}(192.00,120.00)(-16.00,0.00)
\bezier{20}(80.00,20.00)(90.00,10.00)(100.00,0.00)
\bezier{20}(80.00,20.00)(70.00,10.00)(60.00,0.00)
\put(80.00,40.00){\line(0,-1){20.00}}
\put(80.00,40.00){\line(0,1){0.00}}
\put(80.00,80.00){\line(0,-1){40.00}}
\put(80.00,120.00){\line(0,-1){40.00}}
\bezier{30}(80.00,80.00)(120.00,40.00)(160.00,0.00)
\bezier{30}(80.00,80.00)(40.00,40.00)(0.00,0.00)
\end{picture}}
}

\def\ZbbbZb{{
\unitlength=.07pt
\begin{picture}(192.00,120.00)(-16.00,0.00)
\put(40.00,0.00){\line(0,1){20.00}}
\put(40.00,40.00){\line(0,-1){40.00}}
\bezier{20}(40.00,40.00)(60.00,20.00)(80.00,0.00)
\put(80.00,40.00){\line(0,1){0.00}}
\put(80.00,120.00){\line(0,-1){40.00}}
\bezier{30}(80.00,80.00)(120.00,40.00)(160.00,0.00)
\bezier{30}(80.00,80.00)(40.00,40.00)(0.00,0.00)
\end{picture}}
}

\def\bZbbbZ{{
\unitlength=0.07pt
\begin{picture}(192.00,120.00)(-16.00,0.00)
\put(120.00,40.00){\line(0,-1){40.00}}
\bezier{20}(120.00,40.00)(100.00,20.00)(80.00,0.00)
\put(80.00,40.00){\line(0,1){0.00}}
\put(80.00,120.00){\line(0,-1){40.00}}
\bezier{30}(80.00,80.00)(120.00,40.00)(160.00,0.00)
\bezier{30}(80.00,80.00)(40.00,40.00)(0.00,0.00)
\end{picture}}
}

\def\ZvvZv{{
\unitlength=.27pt
\begin{picture}(48.00,30.00)(-4,-30.00)
\bezier{30}(20.00,-20.00)(30.00,-10.00)(40.00,0.00)
\bezier{30}(20.00,-20.00)(10.00,-10.00)(0.00,0.00)
\bezier{20}(10.00,-10.00)(15.00,-5.00)(20.00,0.00)
\put(20.00,-30.00){\line(0,1){10.00}}
\end{picture}}}

\def\coJac#1#2#3{{
\unitlength=.4pt
\begin{picture}(48.00,30.00)(-4,-30.00)
\bezier{30}(20.00,-20.00)(30.00,-10.00)(40.00,0.00)
\bezier{30}(20.00,-20.00)(10.00,-10.00)(0.00,0.00)
\bezier{20}(10.00,-10.00)(15.00,-5.00)(20.00,0.00)
\put(20.00,-30.00){\line(0,1){10.00}}
\put(0,5){\makebox(0,0)[b]{\scriptsize $#1$}}
\put(20,5){\makebox(0,0)[b]{\scriptsize $#2$}}
\put(40,5){\makebox(0,0)[b]{\scriptsize $#3$}}
\end{picture}}}

\def\vZvvZ{{
\unitlength=.27pt
\begin{picture}(48.00,30.00)(-4,-30.00)
\bezier{34}(20.00,-20.00)(30.00,-10.00)(40.00,0.00)
\bezier{34}(20.00,-20.00)(10.00,-10.00)(0.00,0.00)
\bezier{20}(30.00,-10.00)(25.00,-5.00)(20.00,0.00)
\put(20.00,-30.00){\line(0,1){10.00}}
\end{picture}}
}

\def\vZvvZdva{{
\unitlength=.2pt
\begin{picture}(48.00,40.00)(-4.00,0.00)
\bezier{10}(30.00,30.00)(25.00,35.00)(20.00,40.00)
\bezier{34}(0.00,0.00)(20.00,20.00)(40.00,40.00)
\bezier{34}(0.00,40.00)(20.00,20.00)(40.00,0.00)
\end{picture}}}

\def\dvabZbbZ{{
\unitlength=.2pt
\begin{picture}(48.00,40.00)(-4.00,-40.00)
\bezier{10}(30.00,-30.00)(25.00,-35.00)(20.00,-40.00)
\bezier{34}(0.00,0.00)(20.00,-20.00)(40.00,-40.00)
\bezier{34}(0.00,-40.00)(20.00,-20.00)(40.00,0.00)
\end{picture}}}

\def\triZbbZb{{
\unitlength=.2pt
\begin{picture}(48.00,40.00)(-4.00,0.00)
\put(20.00,40.00){\line(0,-1){20.00}}
\bezier{10}(10.00,10.00)(15.00,5.00)(20.00,0.00)
\bezier{34}(0.00,0.00)(20.00,20.00)(40.00,40.00)
\bezier{34}(0.00,40.00)(20.00,20.00)(40.00,0.00)
\end{picture}}}

\def\ZvvZvdva{{
\unitlength=.2pt
\begin{picture}(48.00,40.00)(-4.00,0.00)
\bezier{10}(10.00,30.00)(15.00,35.00)(20.00,40.00)
\bezier{34}(0.00,0.00)(20.00,20.00)(40.00,40.00)
\bezier{34}(0.00,40.00)(20.00,20.00)(40.00,0.00)
\end{picture}}}

\def\tribZbbZ{{
\unitlength=.2pt
\begin{picture}(48.00,40.00)(-4.00,0.00)
\put(20.00,40.00){\line(0,-1){20.00}}
\bezier{10}(30.00,10.00)(25.00,5.00)(20.00,0.00)
\bezier{34}(0.00,0.00)(20.00,20.00)(40.00,40.00)
\bezier{34}(0.00,40.00)(20.00,20.00)(40.00,0.00)
\end{picture}}
}

\def\ZvvZvtri{{
\unitlength=.2pt
\begin{picture}(48.00,40.00)(-4.00,0.00)
\put(20.00,0.00){\line(0,1){20.00}}
\bezier{10}(10.00,30.00)(15.00,35.00)(20.00,40.00)
\bezier{34}(0.00,0.00)(20.00,20.00)(40.00,40.00)
\bezier{34}(0.00,40.00)(20.00,20.00)(40.00,0.00)
\end{picture}}
}

\def\vZvvZtri{{
\unitlength=.2pt
\begin{picture}(48.00,40.00)(-4.00,0.00)
\put(20.00,0.00){\line(0,1){20.00}}
\bezier{10}(30.00,30.00)(25.00,35.00)(20.00,40.00)
\bezier{34}(0.00,0.00)(20.00,20.00)(40.00,40.00)
\bezier{34}(0.00,40.00)(20.00,20.00)(40.00,0.00)
\end{picture}}
}

\def\dvaZbbZb{{
\unitlength=.2pt
\begin{picture}(48.00,40.00)(-4.00,-40.00)
\bezier{10}(10.00,-30.00)(15.00,-35.00)(20.00,-40.00)
\bezier{34}(0.00,0.00)(20.00,-20.00)(40.00,-40.00)
\bezier{34}(0.00,-40.00)(20.00,-20.00)(40.00,0.00)
\end{picture}}}

\def\ixx{\unitlength 1mm 
\linethickness{0.4pt}
\begin{picture}(7,5)(6.5,9)
\put(10,10){\line(0,1){0}}
\multiput(10,10)(.0337079,.0337079){89}{\line(0,1){.0337079}}
\put(13,13){\line(0,1){0}}
\multiput(10,10)(-.0337079,.0337079){89}{\line(0,1){.0337079}}
\multiput(11,13)(-.033333,-.1){30}{\line(0,-1){.1}}
\put(10,10){\line(0,1){0}}
\multiput(10,10)(-.033333,.1){30}{\line(0,1){.1}}
\put(10,10){\line(0,-1){3}}
\put(10,10){\makebox(0,0){$\yyy$}}
\put(12,10){\makebox(0,0)[l]{$x$}}
\end{picture}
}

\def\icc{\unitlength 1mm 
\linethickness{0.4pt}
\begin{picture}(7,5)(6.5,9)
\put(10,10){\line(0,1){0}}
\multiput(10,10)(.0337079,.0337079){89}{\line(0,1){.0337079}}
\put(13,13){\line(0,1){0}}
\multiput(10,10)(-.0337079,.0337079){89}{\line(0,1){.0337079}}
\multiput(11,13)(-.033333,-.1){30}{\line(0,-1){.1}}
\put(10,10){\line(0,1){0}}
\multiput(10,10)(-.033333,.1){30}{\line(0,1){.1}}
\put(10,10){\line(0,-1){3}}
\put(10,10){\makebox(0,0){$\yyy$}}
\put(12,10){\makebox(0,0)[l]{$c$}}
\end{picture}
}
\def\xab{
\unitlength 1mm 
\linethickness{0.4pt}
\begin{picture}(10,5)(9,12.25)
\multiput(10,19)(.0337079,-.0337079){89}{\line(0,-1){.0337079}}
\multiput(13,16)(.0337079,.0337079){89}{\line(0,1){.0337079}}
\multiput(14,19)(-.033333,-.1){30}{\line(0,-1){.1}}
\multiput(13,16)(-.033333,.1){30}{\line(0,1){.1}}
\put(12,19){\line(0,1){0}}
\put(13,10){\line(0,-1){3}}
\bezier{229}(13,16)(17,13)(13,10)
\bezier{229}(13,16)(9,13)(13,10)
\put(13,16){\makebox(0,0){$\yyy$}}
\put(13,10){\makebox(0,0){$\yyy$}}
\put(16,16){\makebox(0,0)[l]{$a$}}
\put(16,10){\makebox(0,0)[l]{$b$}}
\end{picture}
}

\def\aaa{
\unitlength 1mm 
\linethickness{0.4pt}
\begin{picture}(7,5)(6.5,9)
\put(10,10){\line(0,1){0}}
\multiput(10,10)(.0337079,.0337079){89}{\line(0,1){.0337079}}
\put(13,13){\line(0,1){0}}
\multiput(10,10)(-.0337079,.0337079){89}{\line(0,1){.0337079}}
\multiput(11,13)(-.033333,-.1){30}{\line(0,-1){.1}}
\put(10,10){\line(0,1){0}}
\multiput(10,10)(-.033333,.1){30}{\line(0,1){.1}}
\put(10,10){\makebox(0,0){$\yyy$}}
\put(12,10){\makebox(0,0)[l]{$a$}}
\multiput(10,10)(-.0333333,-.05){60}{\line(0,-1){.05}}
\multiput(10,10)(.0333333,-.05){60}{\line(0,-1){.05}}
\end{picture}
}

\def\uuu{
\unitlength 1mm 
\linethickness{0.4pt}
\begin{picture}(7,5)(6.5,9)
\put(10,10){\line(0,-1){3}}
\put(10,10){\makebox(0,0){$\yyy$}}
\put(12,10){\makebox(0,0)[l]{$u$}}
\multiput(10,10)(.0333333,.0666667){60}{\line(0,1){.0666667}}
\multiput(10,10)(-.0333333,.0666667){60}{\line(0,1){.0666667}}
\end{picture}
}

\def\uub{
\unitlength 1mm 
\linethickness{0.4pt}
\begin{picture}(7,5)(6.5,12.25)
\put(12,19){\line(0,1){0}}
\put(13,10){\line(0,-1){3}}
\put(13,10){\makebox(0,0){$\yyy$}}
\put(16,10){\makebox(0,0)[l]{$b$}}
\bezier{170}(13,10)(16,13)(15,16)
\bezier{170}(13,10)(10,13)(11,16)
\multiput(16,19)(-.033333,-.1){30}{\line(0,-1){.1}}
\multiput(15,16)(-.033333,.1){30}{\line(0,1){.1}}
\multiput(12,19)(-.033333,-.1){30}{\line(0,-1){.1}}
\multiput(11,16)(-.033333,.1){30}{\line(0,1){.1}}
\put(10,19){\line(0,1){0}}
\put(11,16){\makebox(0,0){$\yyy$}}
\put(15,16){\makebox(0,0){$\yyy$}}
\put(9.5,16){\makebox(0,0)[r]{$u$}}
\put(17,16){\makebox(0,0)[l]{$u$}}
\end{picture}
}

\newcommand{\abs}[1]{\vert #1 \vert}
\newcommand{\Aut}{{\it Aut}}
\newcommand{\Bij}{{\it Bij}}
\newcommand{\Gr}{{\it Gr}}
\newcommand{\In}{{\it In}}
\newcommand{\Mor}{{\it Mor}}
\DeclareSymbolFont{AMSb}{U}{msb}{m}{n}
\DeclareMathSymbol{\nr}{\mathalpha}{AMSb}{"52}

\newcommand{\Out}{{\it Out}}
\newcommand{\rank}{{\it rank}}
\newcommand{\RG}{{\it RG}}
\newcommand{\rGr}{\overline{\it Gr}}
\newcommand{\rTr}{\overline{\it Tr}}

\title{PROPped up graph cohomology}

\author{M.~Markl\thanks{Partially supported by the grant GA \v CR
                  201/02/1390 and by
   the Academy of Sciences of the Czech Republic,
   Institutional Research Plan No.~AV0Z10190503.}\inst{1}\and
A.A. Voronov\thanks{Partially supported by NSF grant
                  DMS-0227974.}\inst{2}}

\institute{Instutute of Mathematics, Academy of Sciences, \v Zitn\'a 25, 115 67
Praha 1, The~Czech Republic,
\texttt{markl@math.cas.cz}
\and School of Mathematics, University of Minnesota, 206 Church
St. S.E., Minneapolis,
MN 55455, USA, \texttt{voronov@math.umn.edu}}

\begin{document}

\maketitle

\begin{center}
{\it Dedicated to Yuri I. Manin on the occasion of his
seventieth birthday}
\end{center}

\begin{abstract}
  We consider graph complexes with a flow and compute their
  cohomology. More specifically, we prove that for a \PROP\ generated
  by a Koszul dioperad, the corresponding graph complex gives a
  minimal model of the \PROP. We also give another proof of the
  existence of a minimal model of the bialgebra \PROP\ from
  \cite{markl:ba}. These results are based on the useful notion of a
  \hPROP\ introduced by Kontsevich in~\cite{kontsevich:message}.
\end{abstract}

\section*{Introduction}

Graph cohomology is a term coined by M.~Kontsevich
\cite{kontsevich:93,kont:feyn} for the cohomology of
complexes spanned by graphs of a certain type with a differential
given by vertex expansions (also known as splittings), i.e., all
possible insertions of an edge in place of a vertex. Depending on the
type of graphs considered, one gets various ``classical'' types of
graph cohomology. One of them is the graph cohomology implicitly
present in the work of M.~Culler and K.~Vogtmann
\cite{culler-vogtmann:IM86}. It is isomorphic to the rational homology of
the ``outer space,'' or equivalently, the rational homology of the
outer automorphism group of a free group. Another type is the
``fatgraph,'' also known as ``ribbon graph,'' cohomology of
R.C.~Penner \cite{penner:JDG88}, which is isomorphic to the rational
homology of the moduli spaces of algebraic curves.

These types of graph cohomology appear to be impossible to compute, at
least at this ancient stage of development of mathematics. For
example, the answer for ribbon graph cohomology is known only in a
``stable'' limit, as the genus goes to infinity, see a recent ``hard''
proof of the Mumford conjecture by I.~Madsen and M.S. Weiss
\cite{madsen-weiss}. No elementary method of computation seems to
work: the graph complex becomes very complicated combinatorially in
higher degrees. Even the apparently much simpler case of tree
cohomology had been quite a tantalizing problem (except for the
associative case, when the computation follows from the
contractibility of the associahedra) until V.~Ginzburg and
M.M. Kapranov \cite{ginzburg-kapranov:DMJ94} attacked it by
developing Koszul duality for operads.

This paper has originated from a project of computing the cohomology
of a large class of graph complexes. The graph complexes under
consideration are ``\PROP{ped} up,'' which means that the
graphs are directed, provided with a flow, and decorated by the
elements of a certain vector space associated to a given \PROP. When
this \PROP\ is $\sfIB$, the one describing infinitesimal bialgebras,
see M.~Aguiar \cite{aguiar}, we get a directed version of the ribbon
graph complex, while the \PROP\ $\sfLieB$ describing Lie bialgebras
gives a directed commutative version of the graph complex. In both
cases, as well as in more general situations of a directed graph
complex associated to a \PROP\ coming from a Koszul dioperad in the
sense of W.L. Gan \cite{gan} and of a similar graph complex with a
differential perturbed in a certain way, we prove that the
corresponding graph complex is acyclic in all degrees but one, see
Corollary~\ref{37}, answering a question of D.~Sullivan in the Lie
case. This answer stands in amazing contrast with anything one may
expect from the nondirected counterparts of graph cohomology, such as
the ones mentioned in the previous paragraphs: just putting a flow on
graphs in a graph complex changes the situation so dramatically!

Another important goal of the paper is to construct free resolutions
and minimal models of certain \PROP{s}, which might be thought of as
Koszul-like, thus generalizing both the papers of Ginzburg-Kapranov
\cite{ginzburg-kapranov:DMJ94} and Gan \cite{gan}, from trees (and
operads and dioperads, respectively) to graphs (and \PROP{s}). This is
the content of Theorem~\ref{45} below. This theory is essential for
understanding the notions of strongly homotopy structures described by
the cobar construction for Koszul dioperads in \cite{gan} and the
resolution of the bialgebra \PROP\ in~\cite{markl:ba}.

We also observe that axioms of some important algebraic structures
over \PROP{s} can be seen as perturbations of axioms of structures
over \hPROP{s}, objects in a way much smaller than \PROP{s} and even
smaller than dioperads, whose definition, suggested by
Kontsevich~\cite{kontsevich:message}, we give in
Section~\ref{definitions}.  For example, we know from~\cite{markl:ba}
that the \PROP\ $\sfB$ describing bialgebras is a perturbation of the
\hPROP\ $\sfhb$ for $\frac12$bialgebras (more precisely, $\sfB$ is a
perturbation of the \PROP\ generated by the \hPROP\ $\sfhb$). Another
important perturbation of $\sfhb$ is the dioperad ${\it IB}$ for
infinitesimal bialgebras and, of course, also the \PROP\ $\sfIB$
generated by this dioperad. In the same vein, the dioperad ${\it
LieB}$ describing Lie bialgebras and the corresponding \PROP\
$\sfLieB$ are perturbations of the \hPROP\ $\sfhlieb$ for $\frac12$Lie
bialgebras introduced in Example~\ref{79}.

As we argued in~\cite{markl:ba}, every minimal model of a \PROP\ or dioperad
which is a perturbation of a \hPROP\ can be expected to be a
perturbation of a minimal model of this \hPROP. There however might
be some unexpected technical difficulties in applying this principle,
such as the convergence problem in the case of the bialgebra \PROP, 
see Section~\ref{mmod}.

The above observation can be employed to give a new proof of Gan's
results on Koszulness of the dioperads describing Lie bialgebras and
infinitesimal bialgebras. First, one proves that the \hPROP{s}
$\sfhb$ and $\sfhlieb$ are Koszul in the sense of
Section~\ref{koszul}, simply repeating Gan's proof in the simpler case
of \hPROP{s}. This means the \hPROP\ cobar constructions on
the quadratic duals of these \hPROP{s} are minimal models
thereof. Then one treats the dioperadic cobar constructions on the
dioperadic quadratic duals of $\itIB$ and $\itLieB$ as perturbations
of the dg dioperads freely generated by the \hPROP\ cobar
constructions and applies our perturbation theory to show that these
dioperadic cobar constructions form minimal models of the
corresponding dioperads, which is equivalent to their Koszulness.

This paper is based on ideas of the paper \cite{markl:ba} by the first
author and an e-mail message \cite{kontsevich:message} from
Kontsevich. The crucial notion of a \hPROP\ (called
in~\cite{kontsevich:message} a small \PROP) and the idea that
generating a \PROP\ out of a \hPROP\ constitutes a polynomial
functor belong to him.

\noindent
{\em Acknowledgment}: We are grateful to Wee Liang Gan, Maxim
Kontsevich, Sergei Merkulov, Jim Stasheff, and Dennis Sullivan for
useful discussions.

\vskip 3mm
\noindent 
{\bf Table of content:} \ref{definitions}.  
                     \PROP{s}, dioperads and \hPROP{s} 
\hfill\break\noindent 
\hphantom{{\bf Table of content:\hskip .5mm}}  \ref{fprops}.
                     Free \PROP{s}
\hfill\break\noindent 
\hphantom{{\bf Table of content:\hskip .5mm}}  \ref{101}.
                     From \hPROP{s} to \PROP{s}
\hfill\break\noindent 
\hphantom{{\bf Table of content:\hskip .5mm}}  \ref{42}. 
                     Quadratic duality and Koszulness for
                     \hPROP{s}
\hfill\break\noindent 
\hphantom{{\bf Table of content:\hskip .5mm}}  \ref{pert:tech}.
                     Perturbation techniques for graph cohomology 
\hfill\break\noindent 
\hphantom{{\bf Table of content:\hskip .5mm}}  \ref{mmod}.
                     Minimal models of \PROP{s}
\hfill\break\noindent 
\hphantom{{\bf Table of content:\hskip .5mm}}  \ref{classical}.
                     Classical graph cohomology

\section{PROP{s}, dioperads and $\frac12$PROP{s}}
\label{definitions}

Let $k$ denote a ground field which will always be assumed of
characteristic zero. This guarantees the complete reducibility of
finite group representations.
A \PROP\ is a collection $\sfP = \{\sfP(m,n)\}$, $m,n \geq
1$, of {\em differential graded\/} (dg) $(\Sigma_m,\Sigma_n)$-bimodules
(left $\Sigma_m$- right $\Sigma_n$-modules such that the left action
commutes with the right one),
together with two types of compositions, {\em horizontal\/}
\[
\ot : \sfP(m_1,n_1) \ot \cdots \ot \sfP(m_s,n_s) \to 
\sfP(m_1+\cdots +m_s,n_1+\cdots +n_s),
\]
defined for all $m_1,\ldots,m_s,n_1,\ldots,n_s > 0$, and {\em vertical\/}
\[
\circ : \sfP(m,n) \ot \sfP(n,k) \to \sfP(m,k).
\]
defined for all $m,n,k > 0$. These compositions respect the dg
structures. One also assumes the existence of 
a {\em unit\/} $\id \in \sfP(1,1)$. 

\PROP{s} should satisfy axioms which could be 
read off from the example of the {\em
endomorphism \PROP\/} $\End_V$ of a vector space $V$, with
$\End_V(m,n)$ the space of linear maps
$\Hom(\otexp Vn,\otexp Vm)$ with $n$ `inputs' and $m$ `outputs,' 
$\id \in \End_V(1,1)$ the identity map,
horizontal composition given by the tensor product of linear maps, and
vertical composition by the ordinary composition of linear
maps. For a precise definition 
see~\cite{maclane:RiceUniv.Studies63,markl:JPAA96}.

Let us denote, for later use, by $\jcircsoft ji : \sfP(m_1,n_1) \ot
\sfP(m_2,n_2) \to \sfP(m_1+m_2-1,n_1+n_2-1)$, $a,b \mapsto a
\jcircsoft ji b$, $1 \leq i \leq n_1,\ 
1 \leq j \leq m_2$, the operation that composes the 
$j$th output of $b$ to the $i$th input of $a$. Formally,
\begin{equation}
\label{soft}
a \jcircsoft ji b := (\id \ot \cdots \ot \id \ot a \ot \id \ot \cdots
\ot \id)\sigma(\id \ot \cdots \ot \id \ot b \ot \id \ot \cdots \ot
\id),
\end{equation}
where $a$ is at the $j$th place, $b$ is at the $i$th place and
$\sigma \in \Sigma_{n_1+m_2-1}$ is the block permutation 
$((12)(45))_{i-1,j-1,m_2-j,n_1-i}$,
see~\cite{gan}, where this operation was in fact denoted $\jcircsoft
ij$, for details. 

It will also be convenient to introduce special notations for $\jcircsoft 1i$
and $\jcircsoft j1$, namely 
$\circ_i := \jcircsoft 1i:  \sfP(m_1,n_1) \ot \sfP(1,l) \to
\sfP(m_1,n_1+l-1)$, $1 \leq i \leq n_1$, which can be defined simply by
\begin{equation}
\label{1}
a \circ_i b: =
a \circ (\id \ot \cdots \ot \id \ot b \ot \id \ot \cdots \ot \id) \ 
\mbox { ($b$ at the $i$-th position),}
\end{equation}
and, similarly, $\jcirc j := \jcircsoft j1  \sfP(k,1) \ot \sfP(m_2,n_2) 
\to \sfP(m_2 + k -1,n_2)$,
$1 \leq j \leq m_2$, which can be expressed as
\begin{equation}
\label{2}
c \jcirc j d := 
(\id \ot \cdots \ot \id \ot c \ot \id \ot \cdots \ot \id) \circ d\
\mbox { ($c$ at the $j$-th position).}
\end{equation}

A general iterated composition in a \PROP\ is described by a `flow
chart,' which is a not necessarily connected 
graph of arbitrary genus, equipped with a
`direction of gravity' or a `flow,' see Section~\ref{fprops} for more
details. \PROP{s} are in general gigantic objects, with $\sfP(m,n)$
infinite dimensional for any $m$ and $n$.  W.L.~Gan~\cite{gan}
introduced dioperads which avoid this combinatorial explosion.
Roughly speaking, a dioperad is a \PROP\ in which only compositions
along contractible graphs are allowed.

This can be formally expressed by saying that a {\em dioperad\/}  
is a collection $D =\{D(m,n)\}$, $m,n \geq 1$, of dg
$(\Sigma_m,\Sigma_n)$-bimodules with compositions
\[
\jcircsoft ji : D(m_1,n_1) \ot
D(m_2,n_2) \to D(m_1+m_2-1,n_1+n_2-1),
\]
$1 \leq i \leq n_1$, $1 \leq j \leq m_2$,
that satisfy the axioms satisfied by 
operations $\jcircsoft ji$, see~(\ref{soft}), in a general \PROP. 
Gan~\cite{gan} observed that some interesting objects, like Lie
bialgebras or infinitesimal bialgebras, can be defined using algebras
over dioperads.

M.~Kontsevich~\cite{kontsevich:message} suggested even more radical
simplification consisting in  considering objects for which 
only $\circ_i$ and $\jcirc j$ compositions and their iterations are
allowed. More precisely, he suggested:

\begin{definition}
\label{def:small}
A {\em \hPROP\/} is a collection $\sfs =\{\sfs(m,n)\}$ of dg
$(\Sigma_m,\Sigma_n)$-bimodules $\sfs(m,n)$ defined for all pairs of
natural numbers except $(m,n) = (1,1)$, together with compositions
\begin{equation}
\label{circi}
\circ_i:   \sfs(m_1,n_1) \ot \sfs(1,l) 
\to \sfs(m_1,n_1+l-1),\ 1 \leq i \leq n_1,
\end{equation}
and
\begin{equation}
\label{jcirc}
\jcirc j : \sfs(k,1) \ot \sfs(m_2,n_2) 
\to \sfs(m_2 + k -1,n_2),\ 1 \leq j \leq m_2,
\end{equation}
that satisfy the axioms satisfied by 
operations $\circ_i$ and $\jcirc j$, see~(\ref{1}),~(\ref{2}), in a
general \PROP. 
\end{definition}

We suggest as an exercise to unwrap the above definition, write the
axioms explicitly, and compare them to the axioms of a dioperad
in~\cite{gan}. 
Observe that \hPROP{s} cannot have units, because
$\sfs(1,1)$ is not there. 
Later we will also use the notation
\begin{equation}
\label{6}
\circ := \circ_1 = \jcirc 1 : \sfs(k,1) \ot \sfs(1,l) 
\to \sfs(k,l),\ k,l \geq 2.
\end{equation}
The category of \hPROP{s} will be denoted $\cathPROP$.

\begin{example}
{\rm
\label{ex1}
Since \hPROP{s} do not have units, their nature
is close to that of {\em pseudo-operads\/}
~\cite[Definition~1.16]{markl-shnider-stasheff:book}, 
which are, roughly,  operads without units, with
axioms defined in terms of $\circ_i$-operations. More precisely, the
category of \hPROP{s} $\sfs$ with $\sfs(m,n) = 0$ for $m \geq 2$, is
isomorphic to the category of pseudo-operads $\calP$ with
$\calP(0)=\calP(1) = 0$. 
This isomorphism defines a faithful imbedding $\iota :
\catOper \mapsto \cathPROP$ from the category $\catOper$ of
pseudo-operads $\calP$ with $\calP(0) = \calP(1) = 0$ to the category of
\hPROP{s}. To simplify the terminology, by `operad' we will, in this
paper, always understand a pseudo-operad in the above sense.
}
\end{example}

\begin{example}
{\rm
\label{ex2}
Given a \PROP\ $\sfP$, there exists the `opposite' \PROP\
$\sfP^\dagger$ with $\sfP^\dagger(m,n) := \sfP(n,m)$, for each $m,n \geq
1$. A similar duality exists also for dioperads and
\hPROP{s}. 
Therefore one may define another faithful imbedding, $\iota^\dagger :
\catOper \mapsto \cathPROP$, by $\iota^\dagger(\calP) := \iota
(\calP)^\dagger$, where $\iota$ was defined in Example~\ref{ex1}. 
The image of this imbedding consists of all
\hPROP{s} $\sfs$ with $\sfs(m,n) = 0$ for all $n \geq 2$.
}
\end{example}

Every \PROP\ defines a dioperad by forgetting all compositions which
are not allowed in a dioperad. In the same vein, each dioperad
defines a \hPROP\ if we forget all compositions not
allowed in Definition~\ref{def:small}. These observations can be
organized into the following diagram of forgetful functors, 
in which $\catDiop$ denotes the
category of dioperads:
\begin{equation}
\label{4}
\catPROP \stackrel{\Box_1}\longrightarrow \catDiop 
\stackrel{\Box_2}\longrightarrow \cathPROP.
\end{equation}

The left adjoints $F_1 :  \catDiop \to \catPROP$ and $F_2 :
\cathPROP \to \catDiop$ exist by general nonsense. In fact, we give, in
Section~\ref{101}, an explicit description of these
functors. Of primary importance for us will be the composition
\begin{equation}
\label{3}
F := F_1 \circ F_2 : \cathPROP \to \catPROP,
\end{equation}
which is clearly the left adjoint to the forgetful functor $\Box :=
\Box_2 \circ \Box_1:
\catPROP \to \cathPROP$. Given a \hPROP\ $\sfs$, $F(\sfs)$ could be
interpreted as the {\em free \PROP\ generated by the \hPROP~$\sfs$\/}.

Recall~\cite{maclane:RiceUniv.Studies63,markl:JPAA96} 
that an {\em algebra\/} over a \PROP\ $\sfP$ is a morphism  $\sfP \to
\End_V$ of \PROP{s}. The adjoints above offer an
elegant way to introduce algebras over \hPROP{s} and
dioperads: an algebra over a \hPROP\ $\sfs$ is simply an algebra over
the \PROP\ $F(\sfs)$ and, similarly, an algebra over a dioperad $D$ is
defined to be an algebra over the \PROP\ $F_1(D)$.

The following important theorem, whose proof we postpone to
Section~\ref{101}, follows from the fact, observed by M.~Kontsevich
in~\cite{kontsevich:message}, that $F$ and $F_2$ are, in a certain
sense, {\em polynomial functors\/}, see~(\ref{F(s)}) and~(\ref{F_2(s)}).

\begin{theorem}
\label{polynomial}
The functors $F :  \cathPROP \to \catPROP$ and $F_2 : \cathPROP \to
\catDiop$ are exact. This means that they commute with homology, that
is, given a differential
graded \hPROP\ $\sfs$, $H_*(F(\sfs))$ is naturally isomorphic to
$F(H_*(\sfs))$. In particular, for any morphism $\alpha : 
\sfs \to \sft$ of dg \hPROP{s}, the diagram of graded \PROP{s}
\begin{center}
\setlength{\unitlength}{.8cm}
\begin{picture}(5,3.6)
\thinlines

\put(0,3){\makebox(0,0){$H_*(F(\sfs))$}}
\put(5,3){\makebox(0,0){$H_*(F(\sft))$}}
\put(0,0){\makebox(0,0){$F(H_*(\sfs))$}}
\put(5,0){\makebox(0,0){$F(H_*(\sfs))$}}

\put(2.5,3.3){\makebox(0,0)[b]{\scriptsize $H_*(F(\alpha))$}}
\put(-.5,1.5){\makebox(0,0)[r]{\scriptsize $\cong$}}
\put(5.5,1.5){\makebox(0,0)[l]{\scriptsize $\cong$}}
\put(2.5,0.3){\makebox(0,0)[b]{\scriptsize $F(H_*(\alpha))$}}

\put(1.2,0){\vector(1,0){2.5}}
\put(1.2,3){\vector(1,0){2.5}}
\put(0,2.5){\vector(0,-1){2}}
\put(5,2.5){\vector(0,-1){2}}
\end{picture}
\end{center}
is commutative.
A similar statement is also true for $F_2$ in place of $F$.
\end{theorem}

Let us emphasize here that we {\em do not know\/} whether functor
$F_1$ is exact or not. As a consequence of Theorem~\ref{polynomial} we
immediately obtain:

\begin{corollary}
\label{8}
A morphism $\alpha : \sfs \to \sft$ of dg \hPROP{s} is a
homology isomorphism if and only if $F(\alpha) : F(\sfs) \to F(\sft)$ is a
homology isomorphism. A similar statement is also true for~$F_2$.
\end{corollary}

Let us finish our catalogue of adjoint functors by the following definitions.
By a {\em bicollection\/} we mean a sequence $E =
\{E(m,n)\}_{m,n \geq 1}$ of differential graded 
$(\Sigma_m,\Sigma_n)$-bimodules such
that $E(1,1) = 0$. Let us
denote by $\catbCol$ the category of
bicollections. Display~(\ref{4}) then can be completed into the
following diagram of obvious forgetful functors:

\begin{center}
\unitlength=.5cm
\begin{picture}(10.00,5.7)(0.00,0.00)
 \put(0,5){\makebox(0.00,0.00){$\catPROP$}}
 \put(5,5){\makebox(0.00,0.00){$\catDiop$}}
 \put(10.4,5){\makebox(0.00,0.00){$\cathPROP$}}
 \put(5,0){\makebox(0.00,0.00){$\catbCol$}}
 \put(1.2,5){{\vector(1,0){2.5}}}
 \put(6.4,5){{\vector(1,0){2.7}}}
 \put(5,4.5){{\vector(0,-1){4}}}
 \put(.3,4.5){{\vector(1,-1){4}}}
 \put(9.3,4.1){{\vector(-1,-1){3.6}}}
 \put(2.5,5.5){\makebox(0.00,0.00){\scriptsize$\Box_1$}}
 \put(7.5,5.5){\makebox(0.00,0.00){\scriptsize$\Box_2$}}
 \put(2,2.5){\makebox(0.00,0.00)[r]{\scriptsize$\Box_{\tt P}$}}
 \put(8.2,2.5){\makebox(0.00,0.00)[l]{\scriptsize$\Box_{\frac12\tt P}$}}
 \put(5.2,3){\makebox(0.00,0.00)[l]{\scriptsize$\Box_{\tt D}$}}
\end{picture}
\end{center}

Denote finally by $\Gamma_{\tt P}: \catbCol \to \catPROP$, 
$\Gamma_{\tt D}: \catbCol \to \catDiop$ and 
$\Gamma_{\frac12\tt P}: \catbCol \to \cathPROP$
the left adjoints of the functors $\Box_{\tt P}$, $\Box_{\tt D}$ and
$\Box_{\frac12\tt P}$, respectively.

\vskip 2mm
\noindent 
{\bf Notation.}
We will use capital calligraphic letters ${\calP}$, ${\calQ}$, etc. to denote
operads,  small sans serif fonts ${\sf s}$, ${\sf t}$, etc.~to
denote \hPROP{s}, capital italic fonts $S$,
$T$, etc.~to denote dioperads and capital sans serif fonts ${\sf S}$,
${\sf T}$, etc.~to denote \PROP{s}.

\section{Free PROP{s}}
\label{fprops}

To deal with free \PROP{s} and resolutions, we need to fix a suitable
notion of a graph. Thus, in this paper a \emph{graph} or an
$(m,n)$-\emph{graph}, $m,n \ge 1$, will mean a directed (i.e., each
edge is equipped with direction) finite graph satisfying the following
conditions:
\begin{enumerate}

 \item the valence $n(v)$ of each vertex $v$ is at least three;

 \item each vertex has at least one outgoing and at least one incoming
 edge;

 \item there are no directed cycles;

 \item there are precisely $m$ outgoing and $n$ incoming \emph{legs},
 by which we mean edges incident to a vertex on one side and having a
 ``free end'' on the other; these legs are called the \emph{outputs}
 and the \emph{inputs}, respectively;

\item the legs are labeled, the inputs by $\{1, \dots, n\}$, the
outputs by $\{1, \dots, m\}$.

\end{enumerate}
Note that graphs considered are not necessarily connected. Graphs with
no vertices are also allowed. Those will be precisely the disjoint
unions $\uparrow \uparrow \dots \uparrow$ of a number of directed
edges. We will always assume the flow to go from bottom to top, when
we sketch graphs.

Let $v(G)$ denote the set of vertices of a graph $G$, $e(G)$ the set
of all edges, and $\Out (v)$ (respectively, $\In (v)$) the set of
outgoing (respectively, incoming) edges of a vertex $v \in v(G)$.
With an $(m,n)$-graph $G$, we will associate a \emph{geometric
realization} $\abs{G}$, a CW complex whose 0-cells are the vertices of
the graph $G$, as well as one extra 0-cell for each leg, and 1-cells
are the edges of the graph. The 1-cells of $\abs{G}$ are attached to
its 0-cells, as given by the graph.  The \emph{genus} $\gns(G)$ of a
graph $G$ is the first Betti number $b_1 (\abs{G}) = \rank H_1
(\abs{G})$ of its geometric realization. This terminology derives from
the theory of modular operads, but is not perfect, e.g., our genus is
not what one usually means by the genus for ribbon graphs which are
discussed in Section~\ref{classical}.

An \emph{isomorphism} between two $(m,n)$-graphs $G_1$ and $G_2$ is a
bijection between the vertices of $G_1$ and $G_2$ and a bijection
between the edges thereof preserving the incidence relation, the edge
directions and fixing each leg. Let $\Aut (G)$ denote the group of
automorphisms of graph $G$. It is a finite group, being a subgroup of
a finite permutation group.

Let $E = \{E (m,n) \; | \; m,n \geq 1, (m,n) \ne (1,1)\}$, be a
bicollection, see Section~\ref{definitions}. A~standard trick allows
us to extend the bicollection $E$ to pairs $(A,B)$ of finite sets:
\[
E(A,B) := \Bij ([m], A) \times_{\Sigma_m} E(m,n) \times_{\Sigma_n}
\Bij (B, [n]),
\]
where \Bij\ denotes the set of bijections, $[k] = \{1, 2, \dots, k\}$,
and $A$ and $B$ are any $m$- and $n$-element sets, respectively. We
will mostly ignore such subtleties as distinguishing finite sets of
the same cardinality in the sequel and hope this will cause no
confusion. The inquisitive reader may look up an example of careful
treatment of such things and what came out of it in
\cite{getzler-kapranov:CompM98}.

For each graph $G$, define a vector space
\[
E(G) := \bigotimes_{v \in v(G)} E(\Out (v), \In (v)).
\]
Note that this is an unordered tensor product (in other words, a
tensor product ``ordered'' by the elements of an index set), which
makes a difference for the sign convention in graded algebra, see
\cite[page~64]{markl-shnider-stasheff:book}. By definition, $E(\uparrow) =
k$. We will refer to an element of $E(G)$ as a
$G$-\emph{monomial}. One may also think of a $G$-monomial as a
\emph{decorated graph}. Finally, let
\[
\Gamma_{\tt P} (E) (m,n) := \bigoplus_{G \in \Gr(m,n)} E(G)_{\Aut (G)}
\]
be the $(m,n)$-space of the free \PROP\ on $E$ for $m,n \ge 1$, where
the summation runs over the set $\Gr (m,n)$ of isomorphism classes of
all $(m,n)$-graphs $G$ and
\[
E(G)_{\Aut (G)} := E(G)/ \span( ge - e \; | \; g \in \Aut(G), e \in
E(G) )
\]
is the space of \emph{coinvariants} of the natural action of the
automorphism group $\Aut (G)$ of the graph $G$ on the vector space
$E(G)$. The appearance of the automorphism group is due to the fact
that the ``right'' definition would involve taking the colimit over
the diagram of all graphs with respect to isomorphisms, see
\cite{getzler-kapranov:CompM98}. The space $\Gamma_{\tt P} (E) (m,n)$
is a $(\Sigma_m, \Sigma_n)$-bimodule via the action by relabeling the
legs. Moreover, the collection $\Gamma_{\tt P}(E) = \{\Gamma_{\tt
  P}(E)(m,n) \; | \; \hbox{$m,n \ge 1$}\}$ carries a natural \PROP\ structure
via disjoint union of decorated graphs as horizontal composition and
grafting the outgoing legs of one decorated graph to the incoming legs
of another one as vertical composition. The unit is given by $\id \in
k = E(\uparrow)$. The \PROP\ $\Gamma_{\tt P}(E)$ is precisely the free
\PROP\ introduced at the end of Section~\ref{definitions}.

\section{From $\frac12$PROP{s} to PROP{s}}
\label{101}

Let us emphasize that in this article a dg free \PROP\ means a dg
\PROP\ whose underlying (non-dg) \PROP\ is freely generated (by a
\hPROP, dioperad, bicollection,~...) in the
category of (non-dg) \PROP{s}. Such \PROP{s} are sometimes also called
{\em quasi-free\/} or {\em almost-free\/} \PROP{s}. 

We are going to describe the structure of the functors $F :
\cathPROP \to \catPROP$ and $F_2 : \cathPROP \to \catDiop$ and prove
that they commute with homology, i.e., prove Theorem~\ref{polynomial}.
It is precisely the sense of Equations (\ref{s_G}), (\ref{F(s)}), and
(\ref{F_2(s)}) in which we say that the functors $F$ and $F_2$ are
\emph{polynomial}.

Let ${\sf s}$ be a dg \hPROP. Then the dg free \PROP\ $F({\sf
s})$ generated by ${\sf s}$ may be described as follows.  We call an
$(m,n)$-graph $G$, see Section~\ref{fprops}, \emph{reduced}, if it has
no internal edge which is either a unique output or unique input edge
of a vertex. It is obvious
that each graph can be modified to a reduced one by
contracting all the edges violating this condition, i.e., the edges
like this:
\[
{\unitlength=.46pt
\begin{picture}(24.00,30.00)(0.00,3.00)
\put(0.00,0.00){\line(1,0){20.00}}
\bezier{20}(10.00,10.00)(15,5)(20.00,0.00)
\bezier{20}(10.00,10.00)(5,5)(0.00,0.00)
\put(10.00,10.00){\line(0,1){10.00}}
\put(0.00,20.00){\line(1,0){20.00}}
\put(0.00,30.00){\line(1,0){20.00}}
\put(0.00,20.00){\line(0,1){10.00}}
\put(20.00,20.00){\line(0,1){10.00}}
\end{picture}} \qquad \& \qquad
{\unitlength=.46pt
\begin{picture}(24.00,30.00)(0.00,3.00)
\put(0.00,0.00){\line(1,0){20.00}}
\bezier{20}(0.00,30.00)(5,25)(10.00,20.00)
\bezier{20}(20.00,30.00)(15,25)(10.00,20.00)
\put(10.00,10.00){\line(0,1){10.00}}
\put(0.00,10.00){\line(1,0){20.00}}
\put(0.00,30.00){\line(1,0){20.00}}
\put(0.00,0.00){\line(0,1){10.00}}
\put(20.00,0.00){\line(0,1){10.00}}
\end{picture}},
\]
where a triangle denotes a graph with at least one vertex and exactly
one leg in the direction pointed by the triangle, and a box denotes a
graph with at least one vertex. For each reduced graph $G$, define a
vector space
\begin{equation}
\label{s_G}
{\sf s}(G) := \bigotimes_{v \in v(G)} {\sf s}(\Out (v), \In (v)).
\end{equation}
We claim that the \PROP\ $F({\sf s})$ is given by
\begin{equation}
\label{F(s)}
F({\sf s}) (m,n) = \bigoplus_{G \in \rGr(m,n)} {\sf s}(G)_{\Aut (G)},
\end{equation}
where the summation runs over the set $\rGr (m,n)$ of isomorphism
classes of all \emph{reduced} $(m,n)$-\emph{graphs} $G$ and ${\sf
s}(G)_{\Aut (G)}$ is the space of coinvariants of the natural action
of the automorphism group $\Aut (G)$ of the graph $G$ on the vector
space ${\sf s}(G)$.  The \PROP\ structure on the whole collection
$\{F({\sf s}) (m,n)\}$ will be given by the action of the permutation
groups by relabeling the legs and the horizontal and vertical
compositions by disjoint union and grafting, respectively. If grafting
creates a nonreduced graph, we will contract the bad edges and use
suitable \hPROP\ compositions to decorate the reduced graph
appropriately.

A unit in the \PROP\ $F({\sf s})$ is given by $\id \in {\sf s}
(\uparrow)$. A differential is defined as follows. Define a
differential on ${\sf s}(G) = \bigotimes_{v \in v(G)} {\sf s}(\Out
(v), \In (v))$ as the standard differential on a tensor product of
complexes. The action of $\Aut(G)$ on ${\sf s}(G)$ respects this
differential and therefore the space ${\sf s}(G)_{\Aut (G)}$ of
coinvariants inherits a differential. Then we take the standard
differential on the direct sum (\ref{F(s)}) of complexes.

\begin{proposition}
\label{free}
The dg \PROP\ $F({\sf s})$ is the dg free \PROP\ generated by a dg
\hPROP\ {\sf s}, as defined in Section~\ref{definitions}.
\end{proposition}

\noindent 
{\it Proof.}
What we need to prove is that this construction delivers a left
adjoint functor for the forgetful functor $\Box : \catPROP \to
\cathPROP$. Let us define two maps
\[
\Mor_{\cathPROPsoft} ({\sf s}, \Box({\sf P})) 
\arrows_for_Sasha \psi\phi
\Mor_{\catPROP} (F({\sf s}), {\sf P}),
\]
which will be inverses of each other. For a morphism $f: {\sf s} \to
\Box({\sf P})$ of \hPROP{s} and a reduced graph decorated by elements
$s_v \in {\sf s}(m,n)$ at each vertex $v$, we can always compose
$f(s_v)$'s in ${\sf P}$ as prescribed by the graph. The associativity
of \PROP\ compositions in ${\sf P}$ ensures the uniqueness of the
result. This way we get a \PROP\ morphism $\phi(f): F({\sf s}) \to
{\sf P}$.

Given a \PROP\ morphism $g: F({\sf s}) \to {\sf P}$, restrict it to
the sub-\hPROP\ ${\sf s}' \subset F({\sf s})$ given by decorated
graphs with a unique vertex, such as
{\unitlength=.1pt
\begin{picture}(96.00,80.00)(-8.00,0.00)
\bezier{20}(40.00,40.00)(40.00,60.00)(40.00,80.00)
\bezier{30}(0.00,0.00)(40.00,40.00)(80.00,80.00)
\bezier{30}(0.00,80.00)(40.00,40.00)(80.00,0.00)
\bezier{20}(40.00,40.00)(48.50,14.00)(53.50,0.00)
\bezier{20}(40.00,40.00)(32.50,16.00)(26.00,0.00)
\end{picture}}.
We define $\psi(g)$ as the resulting morphism of \hPROP{s}. \qed

\begin{remark}
  {\rm The above construction of the dg free \PROP\ $F({\sf s})$
  generated by a \hPROP\ {\sf s} does not go through for the
  free \PROP\ $F_1 (D)$ generated by a dioperad $D$. The reason is
  that there is no unique way to reduce an $(m,n)$-graph to a graph
  with all possible dioperadic compositions, i.e., all interior edges,
  contracted, as the following figure illustrates:
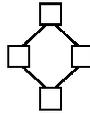
\[
{\unitlength=.8pt
\begin{picture}(44.00,50.00)(0.00,3.00)
\put(15.00,00.00){\line(1,0){10.00}}
\put(15.00,10.00){\line(1,0){10.00}}
\put(0.00,20.00){\line(1,0){10.00}}
\put(0.00,30.00){\line(1,0){10.00}}
\put(30.00,20.00){\line(1,0){10.00}}
\put(30.00,30.00){\line(1,0){10.00}}
\put(15.00,40.00){\line(1,0){10.00}}
\put(15.00,50.00){\line(1,0){10.00}}

\put(0.00,20.00){\line(0,1){10.00}}
\put(10.00,20.00){\line(0,1){10.00}}
\put(15.00,0.00){\line(0,1){10.00}}
\put(25.00,0.00){\line(0,1){10.00}}
\put(30.00,20.00){\line(0,1){10.00}}
\put(40.00,20.00){\line(0,1){10.00}}
\put(15.00,40.00){\line(0,1){10.00}}
\put(25.00,40.00){\line(0,1){10.00}}

\qbezier(7.00,20.00)(12.5,15)(18.00,10.00)
\qbezier(22.00,10.00)(27.5,15)(33.00,20.00)
\qbezier(7.00,30.00)(12.5,35)(18.00,40.00)
\qbezier(22.00,40.00)(27.5,35)(33.00,30.00)
\end{picture}}
\]
This suggests that the functor $F_1$ may be not polynomial.}
\end{remark}

There is a similar description of the dg free dioperad $F_2 ({\sf s})$
generated by a dg \hPROP~{\sf s}:
\begin{equation}
\label{F_2(s)}
F_2 ({\sf s}) (m,n) = \bigoplus_{T \in \rTr(m,n)} {\sf s}(T).
\end{equation}
Here the summation runs over the set $\rTr (m,n)$ of isomorphism
classes of all \emph{reduced} contractible $(m,n)$-graphs $T$. The
automorphism groups of these graphs are trivial and therefore do not show up
in the formula. The following proposition is proven by an obvious
modification of the proof of Proposition~\ref{free}.

\begin{proposition}
  The dg \PROP\ $F_2 ({\sf s})$ is the dg free dioperad generated by a
  dg \hPROP\ {\sf s}, as defined in Section~\ref{definitions}.
\end{proposition}

\noindent
{\it Proof of Theorem \ref{polynomial}}. Let us prove
Theorem~\ref{polynomial} for the dg free \PROP\ $F({\sf s})$ generated
by a dg \hPROP\ {\sf s}. The proof of the statement for $F_2
({\sf s})$ will be analogous and even simpler, because of the absence
of the automorphism groups of graphs.

Proposition~\ref{free} describes $F({\sf s})$ as a direct sum
(\ref{F(s)}) of complexes ${\sf s}(G)_{\Aut (G)}$. Thus the homology
$H_* (F({\sf s}))$ is naturally isomorphic to
\[
\bigoplus_{G \in \rGr(m,n)} H_*( {\sf s}(G)_{\Aut (G)}).
\]
The automorphism group $\Aut(G)$ is finite, acts on ${\sf s}(G)$
respecting the differential, and, therefore, by Maschke's theorem
(remember, we work over a field of characteristic zero), the
coinvariants commute naturally with homology:
\[
H_*( {\sf s}(G)_{\Aut (G)}) \stackrel{\sim}{\to} H_*( {\sf
  s}(G))_{\Aut (G)}.
\]
Then, using the K\"unneth formula, we get a natural isomorphism
\[
H_*( {\sf s}(G)) \stackrel{\sim}{\to} \bigotimes_{v \in v(G)} H_*({\sf
  s}(\Out (v), \In (v))).
\]
Finally, combination of these isomorphisms results in a natural
isomorphism
\[
H_* (F({\sf s)}) \stackrel{\sim}{\to} \bigoplus_{G \in \rGr(m,n)}
\bigotimes_{v \in v(G)} H_*({\sf s}(\Out (v), \In (v)))_{\Aut (G)}
=  F(H_* ({\sf s})).
\]
The diagram in Theorem~\ref{polynomial} is commutative, because of the
naturality of the above isomorphisms.
\qed

\section{Quadratic duality and Koszulness for $\frac12$PROP{s}}
\label{42}\label{koszul}

W.L.~Gan defined in~\cite{gan}, for each dioperad $D$, a dg dioperad
$\Omega_{\tt D}(D) = (\Omega_{\tt D}(D),\partial)$, the {\em cobar
dual\/} of $D$ (${\bf D}D$ in his notation). He also introduced quadratic
dioperads, quadratic duality $D \mapsto D^!$, and showed that, 
for each quadratic
dioperad, there existed a natural map of dg dioperads 
$\alpha_{\tt D}: \Omega_{\tt
D}(D^!) \to D$. He called $D$ {\em Koszul\/}, if $\alpha_{\tt D}$ was a
homology isomorphism. His theory is a dioperadic analog of a similar
theory for operads developed in 1994 by V.~Ginzburg and 
M.M.~Kapranov~\cite{ginzburg-kapranov:DMJ94}. The aim of this section
is to build an analogous theory for \hPROP{s}. Since the passage from
\hPROP{s} to \PROP{s} is given by an exact functor, resolutions of
\hPROP{s} constructed with the help of this theory will induce
resolutions in the category of \PROP{s}.

Let us pause a little and recall, following~\cite{gan}, 
some facts about quadratic duality for dioperads in more
detail. First, a {\em quadratic
dioperad\/} is a dioperad $D$ of the form
\begin{equation}
\label{v_Koline}
D =\Gamma_{\tt D}(U,V)/(A,B,C),
\end{equation}
where $U = \{U(m,n)\}$ is a bicollection with $U(m,n) = 0$ for $(m,n)
\not= (1,2)$, $V = \{V(m,n)\}$ is a bicollection with $V(m,n) = 0$ for
$(m,n) \not= (2,1)$, and $(A,B,C) \subset \Gamma_{\tt D}(U,V)$ denotes
the dioperadic ideal generated by $(\Sigma,\Sigma)$-%
invariant subspaces $A \subset \Gamma_{\tt D}(U,V)(1,3)$, $B \subset
\Gamma_{\tt D}(U,V)(2,2)$ and $C \subset \Gamma_{\tt D}(U,V)(3,1)$.
Notice that we use the original terminology
of~\cite{ginzburg-kapranov:DMJ94} where quadraticity refers to arities
of generators and relations, rather than just relations.  The {\em
  dioperadic quadratic dual\/} $D^!$ is then defined as
\begin{equation}
\label{eLi}
D^! := \Gamma_{\tt D}(U^\vee,V^\vee)/(A^\perp,B^\perp,C^\perp),
\end{equation}
where $U^\vee$ and $V^\vee$ are the linear duals with the action twisted
by the sign representations (the {\em Czech duals\/}, 
see~\cite[p.~142]{markl-shnider-stasheff:book}) and $A^\perp$,
$B^\perp$ and $C^\perp$ are the annihilators of spaces $A$, $B$ and $C$
in 
\[
\Gamma_{\tt D}(U^\vee,V^\vee)(i,j) \cong \Gamma_{\tt
D}(U,V)(i,j)^*,
\] 
where $(i,j) = (1,3)$, $(2,2)$ and $(3,1)$,
respectively. See~\cite[Section~2]{gan} for details.

Quadratic \hPROP{s} and their quadratic duals can then be defined
in exactly the same way as sketched above for dioperads, 
only replacing everywhere $\Gamma_{\tt D}$ by $\Gamma_{\frac12\tt P}$. 
We say therefore that a \hPROP\ $\sfs$ is {\em quadratic\/} if it is
of the form
\[
\sfs =\Gamma_{\frac12\tt P}(U,V)/(A,B,C),
\]
with $U$, $V$ and $(A,B,C) \subset \Gamma_{\frac12\tt P}(U,V)$ 
having a similar obvious meaning as for
dioperads. The {\em quadratic dual\/} of $\sfs$ is defined by a formula completely
analogous to~(\ref{eLi}):
\[
\sfs^! := \Gamma_{\frac12\tt P}(U^\vee,V^\vee)/(A^\perp,B^\perp,C^\perp).
\]
The apparent similarity of the above definitions however hides one very
important subtlety. While
\[
\Gamma_{\tt D}(U^\vee,V^\vee)(1,3) \cong 
\Gamma_{\frac12\tt P}(U^\vee,V^\vee)(1,3)
\]
and
\[
\Gamma_{\tt D}(U^\vee,V^\vee)(3,1) \cong 
\Gamma_{\frac12\tt P}(U^\vee,V^\vee)(3,1),
\]
the $(\Sigma_2,\Sigma_2)$-bimodules
$\Gamma_{\tt D}(U^\vee,V^\vee)(2,2)$ and 
$\Gamma_{\frac12\tt P}(U^\vee,V^\vee)(2,2)$ 
are substantially different, namely
\[
\Gamma_{\tt D}(U^\vee,V^\vee)(2,2) \cong
\Gamma_{\frac12\tt P}(U^\vee,V^\vee)(2,2) \oplus {\rm Ind}^{\Sigma_2 \times
\Sigma_2}_{\{1\}} (U^\vee \ot V^\vee), 
\]
where $\Gamma_{\frac12\tt P}(U^\vee,V^\vee)(2,2) \cong V^\vee \otimes U^\vee$,
see~\cite[section~2.4]{gan} for details. 

We see that the
annihilator of $B \subset \Gamma_{\frac12\tt P}(E,F)(2,2)$ in 
$\Gamma_{\frac12\tt P}(E^\vee,F^\vee)(2,2)$ is much smaller 
than the annihilator of the same
space taken in $\Gamma_{\tt D}(E^\vee,F^\vee)(2,2)$. A consequence of this
observation is a rather stunning fact that quadratic duals {\em do not
commute\/} with functor $F_2 : \cathPROP \to \catDiop$, that is,
$F_2(\sfs^!) \not= F_2(\sfs)^!$. The relation between
$\sfs^!$ and $F_2(\sfs)^!$ is much finer and can be described as follows.

For a \hPROP\ $\sft$, let $\jmath(\sft)$ denote the dioperad which coincides
with $\sft$ as a bicollection and whose structure operations are
those of $\sft$ if they are allowed for \hPROP{s}, and 
are trivial if they are not allowed for \hPROP{s}. This clearly
defines a functor $\jmath : \cathPROP \to \catDiop$.

\begin{lemma}
\label{Eli}
Let $\sfs$ be a quadratic \hPROP. 
Then $F_2(\sfs)$ is a quadratic dioperad and 
\[
F_2(\sfs)^! \cong
\jmath(\sfs^!).
\]
\end{lemma}

\noindent 
{\it Proof.} The proof immediately follows from definitions and we may
safely leave it to the reader.%
\qed

\begin{remark}
{\rm
Obviously $\jmath(\sfs) = F_2(\sfs^!)^!$. This means that the
restriction of the functor $\jmath :
\cathPROP \to \catDiop$ to the full subcategory of quadratic \hPROP{s}
can in fact be defined using quadratic duality. 
}
\end{remark}

The cobar dual $\Omega_{\frac12\tt P}(\sfs)$ of a \hPROP\ $\sfs$ and the
canonical map $\alpha_{\frac12\tt P} : 
\Omega_{\frac12\tt P}(\sfs^!) \to \sfs$ can
be defined by mimicking mechanically the analogous definitions for
dioperads in~\cite{gan}, and we leave this
task to the reader. The following lemma, whose proof is completely
straightforward and hides no surprises, may in fact be interpreted
as a characterization of these objects.

\begin{lemma}
\label{eli1}
For an arbitrary \hPROP\ $\sft$, there exists a functorial isomorphism
of dg dioperads
\[
\Omega_{\tt D}(\jmath(\sft)) \cong 
F_2(\Omega_{\frac12\tt P}(\sft)).
\] 
If $\sfs$ is a quadratic \hPROP, then the canonical maps 
\[
\alpha_{\frac12\tt P} : \Omega_{\frac12\tt P}(\sfs^!) \to \sfs
\] 
and 
\[
\alpha_{\tt D} : \Omega_{\tt D}(F_2(\sfs)^!) \to F_2(\sfs)
\] 
are related by
\begin{equation}
\label{elI}
\alpha_{\tt D} = F_2(\alpha_{\frac12\tt P}).
\end{equation}
\end{lemma}

We say that a quadratic \hPROP\ $\sfs$ is {\em Koszul\/} if the
canonical map $\alpha_{\frac12\tt P} : 
\Omega_{\frac12\tt P}(\sfs^!) \to \sfs$ is
a homology isomorphism. The following proposition is not unexpected,
though it is in fact based on a rather deep Theorem~\ref{polynomial}.

\begin{proposition}
\label{ElI}
A quadratic \hPROP\ $\sfs$ is Koszul if and only if $F_2(\sfs)$ is a
Koszul dioperad.
\end{proposition}

\noindent 
{\it Proof.} The proposition immediately follows from~(\ref{elI}) 
of Lemma~\ref{eli1}
and Corollary~\ref{8} of Theorem~\ref{polynomial}.%
\qed 

We close this section with a couple of important constructions and examples.
Let ${\calP}$ and $\calQ$ be two operads. Recall from Examples~\ref{ex1}
and~\ref{ex2} that ${\calP}$ and $\calQ$ can be considered as \hPROP{s}, 
via imbeddings
$\iota : \catOper \to \cathPROP$ and $\iota^\dagger : \catOper \to
\cathPROP$, respectively. Let us denote
\[
{\calP} * \calQ^\dagger := \iota ({\calP}) \sqcup \iota^\dagger (\calQ)
\]
the coproduct (``free product'') of 
\hPROP{s} $\iota({\calP})$ and $\iota^\dagger(\calQ)$. 
We will need also the quotient
\[
{\calP} \diamond \calQ^\dagger : = 
(\iota({\calP}) \sqcup \iota^\dagger(\calQ))/
(\iota^\dagger(\calQ) \circ \iota({\calP})),
\]
with $(\iota^\dagger(\calQ) \circ \iota({\calP}))$ 
denoting the ideal generated by all
$q^\dagger \circ p$, 
$p \in \iota({\calP})$ and $q^\dagger \in \iota^\dagger(\calQ)$;
here $\circ$ is as in~(\ref{6}).

\begin{exercise}
\label{posledni_den}
{\rm\
Let ${\calP} = \Gamma_{\tt Op}(F)/(R)$ and $\calQ = \Gamma_{\tt Op}(G)/(S)$ be
quadratic operads~\cite[Definition~3.31]{markl-shnider-stasheff:book}; 
here $\Gamma_{\tt Op}(-)$ denotes the free operad
functor. If we interpret $F$, $G$, $R$ and $S$ as bicollections with
\[
F(1,2) := F(2),\ G(2,1) := G(2),\ R(1,3):= R(3) 
\mbox { and } S(3,1) := S(3),
\]
then we clearly have presentations (see~(\ref{v_Koline}))
\[
{\calP}*  \calQ^\dagger = \Gamma_{\frac12\tt P}(F,G)/(R,0,S) 
\mbox { and }
{\calP} \diamond \calQ^\dagger = \Gamma_{\frac12\tt P}(F,G)/(R,G \circ F,S).
\]
which show that both 
${\calP} * \calQ^\dagger$ and ${\calP} \diamond \calQ^\dagger$ are
quadratic \hPROP{s}.
}
\end{exercise}

\begin{exercise}
\label{eli}
{\rm Let $\Ass$ be the operad for associative
algebras~\cite[Definition~1.12]{markl-shnider-stasheff:book}.  Verify
that algebras over \hPROP\ $\Ass* \Ass^\dagger$ are given by a vector
space $V$, an associative multiplication $\cdotvetsi : V \ot V \to V$
and a coassociative comultiplication $\Delta : V \to V \ot V$, with no
relation between these two operations. Verify also that algebras over
$\sfhb := \Ass\diamond \Ass^\dagger$ consists of an associative
multiplication $\cdotvetsi$ and a coassociative comultiplication
$\Delta$ as above, with the exchange rule
\[
\Delta (a\cdotvetsi b) = 0,\ \mbox {for each $a,b \in V$.}
\]
These are exactly {\em $\frac12$bialgebras\/} introduced in~\cite{markl:ba}.
\PROP\ $F(\sfhb)$ generated by \hPROP\ $\sfhb$ is
precisely \PROP\ $\sfhB$ for $\frac 12$bialgebras considered in
the same paper.
}
\end{exercise}

\begin{exercise}
\label{7}
{\rm
Let ${\calP}$ and $\calQ$ be quadratic 
operads~\cite[Definition~3.31]{markl-shnider-stasheff:book}, 
with quadratic duals ${\calP}^!$ and $\calQ^!$, respectively.
Prove that the quadratic dual of the
\hPROP\ ${\calP} \diamond \calQ^\dagger$ is given by 
\[
({\calP} \diamond \calQ^\dagger)^! = {\calP}^! * (\calQ^!)^\dagger.
\]
}
\end{exercise}

\begin{example}
\label{71}
{\rm
The quadratic dual of \hPROP\ $\sfhb$ introduced in
Exercise~\ref{eli} is $\Associative*\Associative^\dagger$. 
Let $\Lie$ denote the operad for 
Lie algebras~\cite[Definition~1.28]{markl-shnider-stasheff:book} and 
$\Com$ the operad
for commutative associative
algebras~\cite[Definition~1.12]{markl-shnider-stasheff:book}.
The quadratic dual of \hPROP\ $\sfhlieb := \Lie \diamond
\Lie^\dagger$ is $\Com * \Com^\dagger$.
}
\end{example}

Gan defined a monoidal structure  $(E,F) \mapsto
E \Box F$ on the category of bicollections such that dioperads were
precisely monoids for this monoidal structure. Roughly speaking, $E
\Box F$ was a sum over all directed contractible graphs $G$ 
equipped with a level function $\ell : v(G) \to \{1,2\}$
such that vertices of level one (that is, vertices with $\ell(v) = 1$)
were decorated by $E$ and vertices of level two were decorated by $F$. 
See~\cite[Section~4]{gan} for precise definitions. 
Needless to say, this $\Box$ should not be mistaken
for the forgetful functors of Section~\ref{definitions}.

Let $D =\Gamma_{\tt D}(U,V)/(A,B,C)$ 
be a quadratic dioperad as in~(\ref{v_Koline}), $\calP :=
\Gamma_{\tt 0p}(U)/(A)$ and $\calQ := \Gamma_{\tt 0p}(V)/(C)$. 
Let us interpret
$\calP$ as a bicollection with $\calP(1,n) = \calP(n)$, $n \geq 1$, 
and let $\calQ^{\rm op}$ be the
bicollection with $\calQ^{\rm op}(n,1) := 
\calQ(n)$, $n \geq 1$, trivial for other
values of $(m,n)$. Since dioperads are $\Box$-monoids in the category of
bicollections, there are canonical maps of bicollections 
\[
\varphi : \calP \Box \calQ^{\rm op} \to D 
\mbox { and } \vartheta :  \calQ^{\rm op} \Box
\calP \to D. 
\]
Let us formulate the following useful proposition.

\begin{proposition}
\label{stavka_pilotu}
The canonical maps 
\[
\varphi : \calP \Box \calQ^{\rm op} \to F_2(\calP \diamond \calQ^\dagger) 
\mbox { and }
\vartheta :  (\calQ^!)^{\rm op} \Box \calP^! 
\to  \calP^! * (\calQ^!)^\dagger 
\]
are isomorphisms of bicollections.
\end{proposition}

\noindent
{\it Proof.}
The fact that $\varphi$ is an isomorphism follows immediately from
definitions. The second isomorphism can be obtained by quadratic
duality: according to~\cite[Proposition~5.9(b)]{gan}, $F_2(\calP \diamond
\calQ^\dagger)^! \cong 
(\calQ^!)^{\rm op} \Box \calP^!$  while $F_2(\calP \diamond
\calQ^\dagger)^! \cong \jmath (\calP^! * 
(\calQ^!)^\dagger) \cong  \calP^! * (\calQ^!)^\dagger$ (isomorphisms
of bicollections) by Lemma~\ref{Eli} and Exercise~\ref{7}.%
\qed

The following theorem is again not surprising, because ${\calP}
\diamond \calQ^\dagger$ 
was constructed from operads ${\calP}$ and $\calQ$ using the relation
\[
q^\dagger \circ p = 0, \mbox { for $p \in {\calP}$ and $q \in \calQ^\dagger$,}
\]
which is a rather trivial {\em mixed distributive law\/} in the sense
of~\cite[Definition~11.1]{fox-markl:ContM97}. 
As such, it cannot create anything unexpected in
the derived category; in particular, it cannot destroy the Koszulness
of ${\calP}$ and $\calQ$. 

\begin{theorem}
\label{eli2}
If ${\calP}$ and $\calQ$ are Koszul quadratic operads, then 
${\calP} \diamond \calQ^\dagger$ is a Koszul \hPROP. 
This implies that the bar construction 
$\Omega_{\frac12\tt P}({\calP}^! * (\calQ^!)^\dagger)$ is a minimal model, in
the sense of~Definition~\ref{mm}, of \hPROP\ ${\calP} \diamond \calQ^\dagger$.
\end{theorem}

\noindent 
{\it Proof.}
We will use the following result of Gan~\cite{gan}.
Given a quadratic dioperad $D$, 
suppose that the operads ${\calP}$ and ${\calQ}$ defined by 
${\calP}(n) := D(1,n)$ and ${\calQ} := D(n,1)$, $n \geq 2$,
are Koszul and that $D \cong {\calP} \Box {\calQ}^{\rm
op}$. Proposition~5.9(c) of~\cite{gan} then states that
$D$ is a Koszul dioperad.

Since, by Proposition~\ref{stavka_pilotu}, 
$F_2({\calP} \diamond \calQ^\dagger) 
\cong {\calP} \Box \calQ^{\rm op}$, 
the above mentioned result implies that
$F_2({\calP} \diamond \calQ^\dagger)$ is a Koszul
quadratic dioperad. Theorem~\ref{eli2} now immediately follows from
Proposition~\ref{ElI} and Exercise~\ref{7}.%
\qed

\begin{example}
\label{Orlik}
{\rm The following example is taken from~\cite{markl:ba}, with signs
altered to match the conventions of the present paper. The minimal
model (see Definition~\ref{mm}) of \hPROP\ $\sfhb$ for
$\frac12$bialgebras, given by the cobar dual 
$\Omega_{\frac12\tt P}(\Associative *
\Associative^\dagger)$, equals
\[
(\Gamma_{\frac12\tt P}(\Xi), \pa_0) 
\stackrel{\alpha_{\frac12\tt P}}{\longrightarrow}
 (\sfhb, \partial = 0) ,
\]
where $\Xi$ denotes the bicollection freely $(\Sigma,\Sigma)$-generated by
the linear span 
$\span(\{\xi^m_n\}_{m,n \in I})$ with
\[
I := \{ m,n \geq 1,\ (m,n) \not= (1,1)\}.
\]
The generator $\xi^m_n$ of biarity $(m,n)$ has degree $n+m-3$. 
The map $\alpha_{\frac12\tt P}$
is defined by
\[
\alpha_{\frac12\tt P}(\xi^1_2) := \jednadva, \
\alpha_{\frac12\tt P}(\xi^2_1) := \dvajedna,
\]
while $\alpha_{\frac12\tt P}$ is trivial on all remaining generators. 
The differential
$\pa_0$ is given by the formula
\begin{eqnarray}
\label{uz_mi_konci_pobyt_v_Minnesote}
\hskip 2mm \pa_0 (\xi^m_n) &:=&
(-1)^{m} \xi_1^m \circ \xi^1_n  +
\sum_U (-1)^{i(s+1) + m + u - s} \xi_u^m \circ_i \xi^1_s 
\\
\nonumber 
&&+ \sum_V (-1)^{(v-j+1)(t+1) - 1} \xi_1^t \hskip 2pt \jcirc j \xi^v_n,
\end{eqnarray}
where we set $\xi^1_1 := 0$,
\[
U := \{u,s \geq 1, \ u+s =n +1, \ 1 \leq i \leq u
\}
\]
and
\[
V = \{t,v \geq 1, \ t+v = m + 1, \ 1 \leq j \leq v
\}.
\]
If we denote $\xi^1_2 = \jednadva$ and $\xi^2_1 = \dvajedna$, then
$\pa_0(\jednadva) = \pa_0(\dvajedna) =0$.
If $\xi^2_2 = \dvadva$, then
\[
\pa_0(\dvadva) = \dvojiteypsilon. 
\]
Under obvious, similar notation,
\begin{eqnarray*}
\pa_0(\jednatri) &=& \ZbbZb - \bZbbZ,
\\
\pa_0(\jednactyri) &=& - \ZbbZbb + \bZbbZb - \bbZbbZ + \ZbbbZb + \bZbbbZ,
\\
\nonumber 
\pa_0(\trijedna) &=& - \ZvvZv + \vZvvZ,
\\
\nonumber 
\pa_0(\dvatri) &=& 
\dvacarkatri - \dvaZbbZb + \dvabZbbZ,
\\
\nonumber 
\pa_0(\tridva) &=&  -
{
\unitlength=.2pt
\begin{picture}(40.00,50.00)(0.00,-50.00)
\put(20.00,-30.00){\line(0,1){10.00}}
\put(20.00,0.00){\line(0,-1){20}}
\bezier{20}(20.00,-20.00)(30.00,-10.00)(40.00,0.00)
\bezier{20}(20.00,-20.00)(10.00,-10.00)(0.00,0.00)
\bezier{20}(20.00,-30.00)(30.00,-40.00)(40.00,-50.00)
\bezier{20}(0.00,-50.00)(10.00,-40.00)(20.00,-30.00)
\end{picture}}
- \ZvvZvdva + \vZvvZdva,
\\
\nonumber 
\pa_0(\tritri)&=& 
-
{
\unitlength=.4pt
\begin{picture}(24.00,30.00)(-2.00,0.00)
\bezier{20}(10.00,10.00)(15.00,5.00)(20.00,0.00)
\bezier{20}(0.00,0.00)(5.00,5.00)(10.00,10.00)
\put(0,-5){
\bezier{20}(10.00,20.00)(15.00,25.00)(20.00,30.00)
\bezier{20}(0.00,30.00)(5.00,25.00)(10.00,20.00)
\put(10.00,30.00){\line(0,-1){25}}
}
\end{picture}}
+ \triZbbZb - \tribZbbZ - \ZvvZvtri + \vZvvZtri,
\\
\nonumber 
\pa_0(\dvactyri) &=& 
{
\unitlength=.1pt
\begin{picture}(96.00,120.00)(-8.00,0.00)
\qbezier(40.00,60.00)(20.00,80.00)(0.00,100.00)
\qbezier(40.00,60.00)(60.00,80.00)(80.00,100.00)
\put(40.00,60.00){\line(0,-1){20.00}}
\qbezier(40.50,39.00)(49.00,14.50)(53.50,0.00)
\qbezier(40.00,40.00)(32.50,16.50)(26.00,0.00)
\qbezier(40.00,40.00)(60.00,20.00)(80.00,0.00)
\qbezier(40.00,40.00)(20.00,20.00)(0.00,0.00)
\end{picture}}
- \dvaZbbbZb - \dvabZbbbZ + \dvaZbbZbb - \dvabZbbZb + \dvabbZbbZ,~\mbox{etc.}
\end{eqnarray*}
}
\end{example}

\begin{example}
\label{79}
{\rm
In this example we discuss a minimal model of the \hPROP\ $\sfhlieb$
introduced in Example~\ref{71}.
The \hPROP\ $\sfhlieb$ describes
{\em $\frac12$Lie bialgebras\/} given by a vector space $V$ with a Lie
multiplication $[-,-] : V \ot V \to V$ and Lie comultiplication
(diagonal)  $\delta : V
\to V \ot V$ tied together by
\[
\delta[a,b] = 0\ \mbox { for all $a,b \in V$.}
\]

A minimal model of $\sfhlieb$ is given 
by the cobar dual $\Omega_{\frac12\tt P}(\Com *\Com^\dagger)$. It is clearly
of the form
\[
(\Gamma_{\frac12\tt P}(\Upsilon), \pa_0) 
\stackrel{\alpha_{\frac12\tt P}}{\longrightarrow}
 (\sfhlieb, \partial = 0) ,
\]
where $\Upsilon$ is the bicollection such that $\Upsilon(m,n)$ is the
ground field placed in degree $m+n-3$ with the sign representation
of $(\Sigma_m,\Sigma_n)$ for $(m,n)\not = 1$, while $\Upsilon(1,1) :=
0$. If we denote by $1^m_n$ the generator of $\Upsilon(m,n)$, then
the map $\alpha_{\frac12\tt P}$ is defined by
\[
\alpha_{\frac12\tt P}(1^1_2) := \jednadva \hskip .5mm , \
\alpha_{\frac12\tt P}(1^2_1) := \dvajedna,
\]
while it is trivial on all remaining generators. There is a formula
for the differential $\pa_0$ which is in fact an anti-symmetric
version of~(\ref{uz_mi_konci_pobyt_v_Minnesote}). We leave writing
this formula, which contains a summation over unshuffles, as an exercise
to the reader.

}
\end{example}

\section{Perturbation techniques for graph cohomology}
\label{pert:tech}

Let $E$ be a bicollection. We are going to introduce, for an arbitrary
fixed $m$ and $n$, three very important gradings of the piece
$\Gamma_{\tt P}(E)(m,n)$ of the free \PROP\ $\Gamma_{\tt P}(E)$. We
know, from Section~\ref{fprops}, that $\Gamma_{\tt P}(E)(m,n)$ is the
direct sum, over the graphs $G \in \Gr(m,n)$, of the vector spaces
$E(G)_{\Aut (G)}$.  Recall that we refer to elements of $E(G)_{\Aut
(G)}$ as {\em $G$-monomials\/}.

The first two gradings are of a purely topological nature.  The {\em
component grading\/} of a $G$-monomial $f$ is defined by $\comp(f) :=
\comp(G)$, where $\comp(G)$ is the number of connected components of
$G$ minus one.  The {\em genus grading\/} is given by the topological
genus $\gns(G)$ of graphs (see Section~\ref{fprops} for a precise
definition), that is, for a $G$-monomial $f$ we put $\gns(f) :=
\gns(G)$. Finally, there is another {\em path grading\/}, denoted
$\pth(G)$, implicitly present in~\cite{kontsevich:message}, defined as
the total number of directed paths connecting inputs with outputs of
$G$.  It induces a grading of $\Gamma_{\tt P}(E)(m,n)$ by setting
$\pth(f) := \pth(G)$ for a $G$-monomial $f$.

\begin{exercise}
\label{133}
{\rm
Prove that for each  $G$-monomial $f \in \Gamma_{\tt P}(E)(m,n)$,
\[
\gns(f) + \max\{m,n\} \leq \pth(f) \leq mn (\gns(f)+1)
\]
and
\[
\comp(f) \leq \min\{m,n\} -1.
\]
Find examples that show that these inequalities cannot be improved and
observe that our assumption that $E(m,n)$ is nonzero only for $m,n
\geq 1$, $(m,n) \not= (1,1)$, is crucial.
}
\end{exercise}

Properties of these gradings are summarized in the following
proposition.

\begin{proposition}
\label{30}
Suppose $E$ is a bicollection of finite dimensional
$(\Sigma,\Sigma)$-bimodules. Then, for any fixed $d$, the subspaces
\begin{equation}
\label{10}
\span\{f \in \Gamma_{\tt P}(E)(m,n); \ \gns(f) = d\}
\end{equation}
and
\begin{equation}
\label{10bis}
\span\{f \in \Gamma_{\tt P}(E)(m,n); \ \pth(f) = d\},
\end{equation}
where $\span\{-\}$ is the $k$-linear span,
are finite dimensional. 
The subspace $\Gamma_{\tt D}(E)(m,n) \subset
\Gamma_{\tt P}(E)(m,n)$ can be characterized as
\begin{equation}
\label{11}
\Gamma_{\tt D}(E)(m,n) = \span\{f \in \Gamma_{\tt P}(E)(m,n);\ 
\comp(f) = \gns(f) = 0\}.
\end{equation}

{}For each $f \in \Gamma_{\tt D}(E)(m,n)$, $\pth(f) \leq mn$, and 
the subspace $\Gamma_{\frac12\tt P}(E)(m,n) \subset
\Gamma_{\tt D}(E)(m,n)$ can be described as
\begin{equation}
\label{12}
\Gamma_{\frac12\tt P}(E)(m,n) 
= \span\{f \in \Gamma_{\tt D}(E)(m,n);\ \pth(f) = mn\}.
\end{equation}
\end{proposition}

\noindent
{\it Proof.}
Since all vertices of our graphs are at least trivalent, it follows
from standard combinatorics that there is only
a finite number of $(m,n)$-graphs with a fixed genus. This proves the
finite-dimensionality of the space in~(\ref{10}). 
Description~(\ref{11}) follows immediately from
the definition of a dioperad.
Our proof of the finite-dimensionality of the space in~(\ref{10bis}) 
is based on the following
argument taken from~\cite{kontsevich:message}.

\begin{figure}[t]
\begin{center}
\unitlength=.6cm
\begin{picture}(5,7.9)(-1,-1)
\thicklines
\put(1,1){\makebox(0.00,0.00){$\bullet$}}
\put(1,3){\makebox(0.00,0.00){$\bullet$}}
\put(2,6){\makebox(0.00,0.00){$\bullet$}}
\put(-1,-1){\line(1,1){2}}
\put(3,-1){\line(-1,1){2}}
\put(1,1){\line(0,1){2}}
\put(1,3){\line(1,1){2}}
\put(1,3){\line(-1,1){1}}
\put(0,4){\vector(1,1){3}}
\put(3,5){\vector(-1,1){2}}
\put(-.5,0){\makebox(0.00,0.00)[r]{$p_1$}}
\put(2.5,0){\makebox(0.00,0.00)[l]{$p_2$}}
\put(.2,1){\makebox(0.00,0.00)[l]{$w$}}
\put(.2,3){\makebox(0.00,0.00)[l]{$v$}}
\put(1.2,6){\makebox(0.00,0.00)[l]{$u$}}
\end{picture}
\end{center}
\caption{\label{fig:ELI}%
Three branching points $u$, $v$ and $w$ of paths $p_1$ and $p_2$.
}
\end{figure}
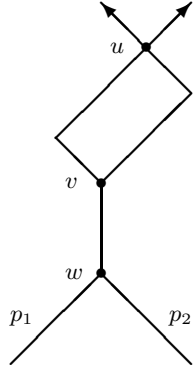

Let us say that a vertex $v$ is a {\em branching vertex\/} for a pair
of directed paths $p_1,p_2$ of a graph $G \in \Gr(m,n)$, if $v$ 
is a vertex of both $p_1$ and $p_2$ and if it has the
property that either there exist two different input edges $f_1,f_2$
of $v$ such that $f_s \in p_s$, $s=1,2$, or there exist two different
output edges $e_1,e_2$ of $v$ such that $e_s \in p_s$, $s=1,2$. 
See also Figure~\ref{fig:ELI}.
Denote
${\br}(p_1,p_2)$ the number of all branching vertices for
$p_1$ and $p_2$. A moment's reflection convinces us 
that a pair of paths $p_1$ and $p_2$ with $b$ branching points
generates at least $2^{b-1}$ different paths in $G$, therefore 
$2^{\br(p_1,p_2)-1} \leq d$, where $d$ is the total number of directed
paths in $G$. This implies that
\[
\br(p_1,p_2) \leq \log_2(d) +1.
\]
Now observe that each vertex is a branching point for at least one
pair of paths. We conclude that the number of vertices
of $G$ must be less than or equal to $d^2 \cdot (\log_2(d)+1)$.

The graph $G$ cannot have
vertices of valence bigger than $d$, because each vertex of valence
$k$ generates at least $k-1$ different paths in $G$.
Since there are only finitely many isomorphism classes of graphs with
the number of vertices bounded by a constant and with the
valences of its vertices bounded by another constant, the
finite-dimensionality 
of the space in~(\ref{10bis}) is proven.

Let us finally demonstrate~(\ref{12}). Observe first that for a graph $G \in
\Gr(m,n)$ of genus $0$, $mn$ is actually an upper bound for
$\pth(G)$, because for each output-input pair $(i,j)$ there exists at most one
path joining $i$ with $j$ (genus $0$ assumption). It is also not
difficult to see that $\pth(f) = mn$ for a $G$-monomial $f \in
\Gamma_{\frac12\tt P}(E)$. So it remains to prove that $\pth(f) = mn$ 
implies $f \in \Gamma_{\frac12\tt P}(E)$.

\begin{figure}[t]
\begin{center}
\unitlength=.6cm
\begin{picture}(5,9.5)(0.00,-5.5)
\thicklines
\put(2.00,00){\makebox(0.00,0.00){$\bullet$}}
\put(4.00,-2){\makebox(0.00,0.00){$\bullet$}}
\put(0.00,-2.00){\vector(1,1){2}} 
\put(4.00,-2.00){\vector(-1,1){2}}
\put(4.00,-2.00){\vector(1,1){2}}
\put(2.00,-0.00){\vector(0,1){2}}
\put(2.00,-0.00){\vector(0,1){2}}
\put(4.00,-4.00){\vector(0,1){2}}
\put(1.70,1){\makebox(0.00,0.00)[r]{$a$}}
\put(3.70,-3){\makebox(0.00,0.00)[r]{$b$}}
\put(.50,-1){\makebox(0.00,0.00)[r]{$e$}}
\put(2.50,-1){\makebox(0.00,0.00)[r]{$f$}}
\put(5.50,-1){\makebox(0.00,0.00)[l]{$g$}}
\put(2.5,0){\makebox(0.00,0.00)[l]{$u$}}
\put(4.5,-2){\makebox(0.00,0.00)[l]{$v$}}
\put(0,-5.5){\makebox(0.00,0.00)[t]{$j_1$}}
\put(4,-5.5){\makebox(0.00,0.00)[t]{$j_2$}}
\put(2,3.5){\makebox(0.00,0.00)[b]{$i_1$}}
\put(6,3.5){\makebox(0.00,0.00)[b]{$i_2$}}
\thinlines
\put(0.00,-2.00){\line(0,-1){3}}
\put(4.00,-2.00){\line(0,-1){3}}
\put(6.00,0.00){\line(0,1){3}}
\put(2.00,0.00){\line(0,1){3}}
\end{picture}
\end{center}
\caption{\label{fig:1}%
A configuration forbidden for \hPROP{s} -- $f$ is a `bad'
edge. Vertices $u$ and $v$
might have more input or output edges which we did not indicate.}
\end{figure}
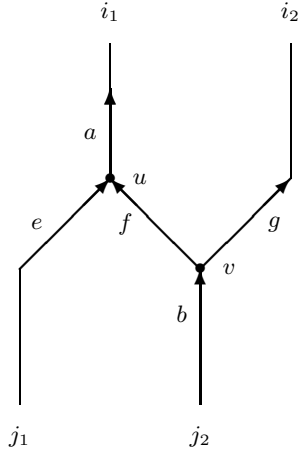

Suppose that $f$ is a $G$-monomial such that 
$f \in \Gamma_{\tt D}(E)(m,n) \setminus \Gamma_{\frac12\tt P}(E)(m,n)$. This
happens exactly when $G$ contains a
configuration shown in Figure~\ref{fig:1}, forbidden for \hPROP{s}.
Then there certainly exists a path $p_1$ containing edges $e$ and $a$,
and another path $p_2$ containing edges $b$ and $g$. Suppose
that $p_s$ connects output $i_s$ with input $j_s$, $i = 1,2$, as in
Figure~\ref{fig:1}. 
It is then clear that
there is no path that connects $i_2$ with $j_1$, which means that the total
number of paths in $G$ is not maximal. This finishes the proof of the
proposition.%
\qed

\begin{remark}
{\rm As we already know, there are various `restricted' versions of
\PROP{s} characterized by types of graphs along which the composition
is allowed.  Thus \hPROP{s} live on contractible 
graphs without `bad' edges as in Figure~\ref{fig:1},
and Gan's dioperads live on all contractible graphs. A version
of \PROP{s} for which only compositions along {\em connected\/} graphs are
allowed was studied by Vallette who called these \PROP{s} {\em
properads\/}~\cite{vallette:thesis}. All this can be summarized by a 
chain of inclusions of full subcategories
\[
\catOper \subset \cathPROP \subset \catDiop \subset \catProper \subset
\catPROP.
\]
}
\end{remark}

Let $\Gamma_\pth(E) \subset \Gamma_{\tt P}(E)$ be the subspace spanned
by all $G$-monomials such that $G$ is contractible and contains at
least one `bad' edge as in Figure~\ref{fig:1}.  By
Proposition~\ref{30}, one might equivalently define $\Gamma_{\pth}(E)$
by
\begin{eqnarray*}
\lefteqn{
\Gamma_{\pth}(E)(m,n)=}
\\ 
&=&
\span\{f \in \Gamma_{\tt D}(E)(m,n);\ \comp(f) =\gns(f) = 0,\
\mbox { and } \pth(f) < mn\}.
\end{eqnarray*}
If we denote
\[
\Gamma_\gencomp(E) := \span\{ f \in \Gamma_{\tt P}(E);\ \comp(f) + \gns(f)
> 0\} 
\] 
then there is a natural decomposition
\[
\Gamma_{\tt P}(E) = \Gamma_{\frac12\tt P}(E) \oplus \Gamma_\pth(E) \oplus
\Gamma_\gencomp(E) 
\]
in which clearly $\Gamma_{\frac12\tt P}(E) \oplus \Gamma_\pth(E) =
\Gamma_{\tt D}(E)$. Let $\pi_{\frac12\tt P}$, 
$\pi_\pth$ and $\pi_\gencomp$ denote the corresponding
projections. For a degree $-1$ differential $\pa$ on
$\Gamma_{\tt P}(E)$, introduce derivations $\pa_0$, $\pa_\pth$ and
$\pa_\gencomp$ determined by their restrictions to the generators $E$ as
follows:
\[
\pa_{0}|_E := \pi_{\frac12\tt P} \circ \pa|_E,\
\pa_{\pth}|_E := \pi_{\pth} \circ \pa|_E \  \mbox { and }
\pa_{\gencomp}|_E := \pi_{\gencomp} \circ \pa|_E.
\] 
Let us define also $\pa_{\tt D}: = \pa_0 + \pa_\pth$, 
the {\em dioperadic part\/} of $\pa$. The decompositions
\begin{equation}
\label{fund}
\pa = \pa_{\tt D} + \pa_{\gencomp} = \pa_{0} + \pa_{\pth}  + \pa_{\gencomp}
\end{equation}
are fundamental for our purposes. We will call them the {\em canonical
decompositions\/} of the differential $\pa$. The following example shows
that, in general, $\pa_0$, $\pa_{\tt D}$ and $\pa_{\gencomp}$ {\em need
not\/} be differentials, as they may not square to zero.

\begin{example}
{\rm Let us consider the free \PROP\ $\Gamma_{\tt P}(a,b,c,u,x)$, where
the generator $a$ has degree $1$ and biarity $(4,2)$, $b$ degree $0$
and biarity $(2,1)$, $c$ degree $1$ and biarity $(4,1)$, $u$ degree
$0$ and biarity $(2,1)$, and $x$ degree $2$ and biarity $(4,1)$.
Define a degree $-1$ differential $\pa$ by the following formulas
whose meaning is, as we believe, clear:
\[
\pa \left( \ixx \right) := \! \xab+ \icc \hskip .1em,\
\pa \left( \aaa \right) := \!\uuu \hskip .3em \otimes \hskip -0.3em \uuu\ ,\
\pa \left( \icc \right) := - \!\uub  \hskip 1.5em,
\]
while $\pa (b) = \pa(u) = 0$. One can easily verify that $\pa^2 =
0$. By definition,
\[
\pa_0 \left( \ixx \right) =  \icc \hskip .5em,\
\pa_0 \left( \aaa \right) = 0 ,\
\pa_0 \left( \icc \right) = - \uub 
\]
and, of course, $\pa_0 (b) = \pa_0(u) = 0$. A simple calculation shows
that
\[
\pa_0^2 \left( \ixx \right) = - \uub
\]
therefore $\pa_0^2 \not= 0$. Since $\pa_0 = \pa_{\tt D}$, we
conclude that also $\pa_{\tt D}^2 \not= 0$.
}
\end{example}

Let us formulate some conditions which guarantee that the
derivations $\pa_0$ and $\pa_{\tt D}$ square to zero.
We say that a differential $\pa$ in $\Gamma_{\tt P}(E)$ is {\em
connected\/} if $\comp(\pa(e))=0$ for each $e \in E$. Similarly we say that
$\pa$ {\em has genus zero\/} if $\gns(\pa(e))=0$ for $e \in E$.
Less formally, connectivity of $\pa$ means that $\pa(e)$ is a sum of
$G$-monomials with all $G$'s connected, and $\pa$ has genus zero
if $\pa(e)$ is a sum of $G$-monomials with all $G$'s of genus $0$ (but
not necessarily connected).

\begin{proposition}
\label{18}
In the canonical decomposition~(\ref{fund}) of a differential $\pa$ in
a free \PROP\ $\Gamma_{\tt P}(E)$, $\pa_{\tt D}^2 = 0$ always
implies that $\pa_0^2 = 0$.

If moreover either (i) the differential $\pa$ is connected or (ii)
$\pa$ has genus zero, then $\pa_{\tt D}^2 = 0$, therefore both $\pa_0$
and $\pa_{\tt D}$ are differentials on $\Gamma_{\tt P}(E)$.
\end{proposition}

\noindent 
{\it Proof.}  For a $G$-monomial $f$, write
\begin{equation}
\label{00}
\pa_{\tt D}(f) = \sum_{H \in U} g_H,
\end{equation}
the sum of $H$-monomials $g_H$ over a finite set $U$ of graphs. Since
$\pa_{\tt D}$ is a derivation, each $H \in U$ is obtained by replacing
a vertex $v \in v(G)$ of biarity $(s,t)$ by a graph $R$ of the same
biarity. It follows from the definition of the dioperadic part
$\pa_{\tt D}$ that each such $R$ is contractible.  This implies that
all graphs $H \in U$ which nontrivially contribute to the
sum~(\ref{00}) have the property that $\pth(H) \leq \pth(G)$ ({\em
$\pa_{\tt D}$ does not increase the path grading\/}) and that
\begin{equation}
\label{33}
\pa_0(f) = \sum_{H \in U_{\tt 0}} g_H, \mbox { where $U_{0} := \{H
\in U;\ \pth(H) = \pth(G)\}$.}
\end{equation}

This can be seen as follows.  It is clear that a replacement of a
vertex by a contractible graph cannot increase the total number of
paths in $G$. This implies that $\pa_{\tt D}$ does not increase the
path grading.  Equation~(\ref{33}) follows from the observation that
decreasing the path grading locally at a vertex decreases the path
grading of the whole graph. By this we mean the following.

Assume that a vertex $v$ of biarity $(s,t)$ is replaced by a
contractible graph $R$ such that $\pth(R) < st$. This means that there
exists an output-input pair $(i,j)$ of $R$ for which there is no path
in $R$ connecting output $i$ with input $j$. On the other hand, in $G$
there certainly existed a path that ran through output $i$ and input
$j$ of vertex $v$ and broke apart when we replaced $v$ by $R$.

Now we see that $\pa_{0}^2$ is precisely the part of $\pa^2_{\tt D}$
which preserves the path grading. This makes the implication
$(\pa_{\tt D}^2 = 0) \Longrightarrow (\pa_0^2 = 0)$ completely obvious and
proves the first part of the proposition.

For the proof of the second half, it will be convenient to
introduce still another grading by putting
\begin{equation}
\label{l}
\grad(G) :=\comp(G)  -\gns(G) = +|v(G)| - |e(G)| -1,
\end{equation}
where $|v(G)|$ denotes the number of vertices and $|e(G)|$ the number
of internal edges of $G$.  Let $f$ be a $G$-monomial as above.  Let us
consider a sum similar to~(\ref{00}), but this time for the entire
differential $\pa$:
\[
\pa(f) = \sum_{H \in S} g_H,
\]
where $S$ is a finite set of graphs. We claim that, under
assumptions~(i) or~(ii),
\begin{equation}
\label{32}
\pa_{\tt D}(f) = \sum_{H \in S_{\tt D}} g_H, \mbox { where $S_{\tt D} := \{H
\in S;\ \grad(H) = \grad(G)\}$.}
\end{equation}
This would clearly imply that $\pa_{\tt D}^2$ is exactly the part of
$\pa^2$ that preserves the $\grad$-grading, therefore $\pa_{\tt D}^2 =
0$.

As in the first half of the proof, each $H \in S$ is obtained from $G$
by replacing $v \in v(G)$ by some graph $R$. In case~(i), $R$ is
connected, that is $\comp(G) = \comp(H)$ for all $H \in S$. It
follows from an elementary algebraic topology that $\gns(H) \geq
\gns(G)$ and that $\gns(G) = \gns(H)$ if and only if $\gns(R) =
0$. This proves~(\ref{32}) for connected differentials.

Assume now that $\pa$ has genus zero, that is, $\gns(R) = 0$. This
means that $R$ can be contracted to a disjoint $R'$ union of
$\comp(R)+1$ corollas. Since $\grad(-)$ is a topological invariant, we
may replace $R$ inside $H$ by its contraction $R'$. We obtain a
graph $H'$ for which $\grad(H) = \grad(H')$. It is obvious
that $H'$ has the same number of internal edges as $G$ and that
$|v(H')| = |v(G)| + \comp(R)$, therefore $\grad(G) = \grad(H) +\comp(R)$.
This means that $\grad(G) = \grad(H)$ if and only if $\comp(R) = 0$, i.e.~if
$R$ is connected. This proves~(\ref{32}) in case~(ii) and finishes the
proof of the Proposition.%
\qed

The following theorem will be our basic tool to calculate the homology
of free differential graded \PROP{s} in terms of the canonical
decomposition of the differential.

\begin{theorem}
\label{35}
Let $(\Gamma_{\tt P}(E),\pa)$ be a dg free \PROP\ and $m,n$ fixed
natural numbers. 

(i) 
Suppose that the differential $\pa$ is connected. 
Then the genus grading defines, by
\begin{equation}
\label{fil1}
F^\gns_p := \span\{f \in \Gamma_{\tt P}(E)(m,n);\  \gns(f) \geq -p\},
\end{equation}
an increasing $\pa$-invariant
filtration of $\Gamma_{\tt P}(E)(m,n)$.

(ii) 
If the differential $\pa$ has genus zero, then 
\[
F^\grad_p := \span\{f \in \Gamma_{\tt P}(E)(m,n);\  \grad(f) \geq -p\}.
\]
is also an increasing $\pa$-invariant filtration of $\Gamma_{\tt
P}(E)(m,n)$.

The spectral sequences induced by these filtrations have both the
first term isomorphic to $(\Gamma_{\tt P}(E)(m,n),\pa_{\tt D})$ and
they both abut to $H_*(\Gamma_{\tt P}(E)(m,n),\pa)$.

(iii)
Suppose that $\pa_{\tt D}^2 = 0$. 
Then the path grading defines an increasing $\pa_{\tt D}$-invariant
filtration 
\[
F^\pth_p := \span\{f \in \Gamma_{\tt P}(E)(m,n);\  \pth(f) \leq p\}.
\]
This filtration induces a first quadrant spectral sequence whose first
term is isomorphic to $(\Gamma_{\tt P}(E)(m,n),\pa_0)$ and which
converges to $H_*(\Gamma_{\tt P}(E)(m,n),\pa_{\tt D})$.
\end{theorem}

\noindent 
{\it Proof.}  The proof easily follows from Proposition~\ref{18} and
the analysis of the canonical decomposition given in the proof of that
proposition.%
\qed

The following proposition describes an important particular case when
the spectral sequence induced by the filtration~(\ref{fil1}) converges.

\begin{proposition}
\label{dnes_prileti_E}
If $\pa$ is connected and preserves the path
grading, then the filtration~(\ref{fil1}) induces a second quadrant
spectral sequence whose first term is isomorphic to $(\Gamma_{\tt
P}(E)(m,n),\pa_{0})$ and which converges to  
$H_*(\Gamma_{\tt P}(E)(m,n),\pa)$.
\end{proposition}

\noindent 
{\it Proof.}  Under the assumptions of the proposition, the path
grading is a $\pa$-invariant grading, compatible with the genus
filtration~(\ref{fil1}), by finite dimensional pieces, see
Proposition~\ref{30}. This guarantees that the generally ill-behaved
second quadrant spectral sequence induced by~(\ref{fil1}) converges.
The proof is finished by observing that the assumption that $\pa$
preserves the path grading implies that $\pa_0 = \pa_{\tt D}$.%
\qed

In most applications either $\pa$ is connected or
$\pa = \pa_{\tt D}$, though there are also natural examples of
\PROP{s} with disconnected differentials, such as the {\em deformation
quantization\/} \PROP\ ${\sf DefQ}$ introduced by Merkulov
in~\cite{merkulov:defquant}.  
The following corollary immediately follows from
Theorem~\ref{35}(iii) and Proposition~\ref{dnes_prileti_E}.

\begin{corollary}
\label{37}
Let $\sfP$ be a graded \PROP\  concentrated in degree $0$ 
and $\alpha : (\Gamma_{\tt P}(E),\pa) \to (\sfP,0)$ a homomorphism of
dg \PROP{s}. Suppose that $\alpha$ induces an isomorphism 
$H_0(\Gamma_{\tt P}(E),\pa) \cong \sfP$ and that  
$\Gamma_{\tt P}(E)$ is $\pa_0$-acyclic in positive
degrees. Suppose moreover that either

(i) 
$\pa$ is connected and preserves the path grading, or

(ii)
$\pa(E) \subset \Gamma_{\tt D}(E)$.

\noindent 
Then $\alpha$ is a free resolution of the \PROP~$\sfP$. 
\end{corollary}

\begin{remark}
{\rm
In Corollary~\ref{37} we assumed that the \PROP\ $\sfP$ was concentrated in
degree~$0$. The case of a general nontrivially graded non-differential
\PROP\ $\sfP$ can be treated by introducing the 
{\em Tate-Jozefiak grading\/}, as it
was done, for example, for bigraded models of operads 
in~\cite[page~1481]{markl:zebrulka}.
}
\end{remark}

\section{Minimal models of PROP{s}}
\label{mmod}

In this section we show how the methods of this paper can be used to 
study minimal models of \PROP{s}. Let us first give a precise
definition of this object.

\begin{definition}
\label{mm}
A {\em minimal model\/} of a dg \PROP\ $\sfP$
is a dg free \PROP\ $(\Gamma_{\tt P}(E),\pa)$
together with a homology isomorphism
\[
\sfP \stackrel{\alpha}{\longleftarrow} (\Gamma_{\tt P}(E),\pa).
\] 
We also assume that the image of $\pa$ consists of 
decomposable elements of $\Gamma_{\tt P}(E)$ or, equivalently, that
$\pa$ has no ``linear part'' (the {\em minimality condition\/}). Minimal models
for \hPROP{s} and dioperads are defined in exactly the same way, only
replacing $\Gamma_{\tt P}(-)$ by 
$\Gamma_{\frac12\tt P}(-)$ or $\Gamma_{\tt D}(-)$.
\end{definition}

The above definition generalizes minimal models for operads 
introduced in~\cite{markl:zebrulka}. While we proved, 
in~\cite[Theorem~2.1]{markl:zebrulka} that each operad admits, under
some very mild conditions, a~minimal model, and while the same
statement is probably true also for dioperads, 
a similar statement for a general \PROP\ would require some way to
handle a divergence problem (see also
the discussion in~\cite{markl:ba} and below).

\vskip 3mm
\noindent 
{\bf Bialgebras.}
Recall that a {\em bialgebra\/} is a
vector space $V$ with an associative multiplication $\cdot : V \ot V
\to V$ and a coassociative comultiplication $\Delta : V \to V \ot V$
which are related by
\[
\Delta(a\cdot b) = \Delta(a) \cdot \Delta (b), \mbox { for } a,b \in V.
\]
The \PROP\ $\sfB$ describing bialgebras has a presentation
$
\sfB = \Gamma_{\tt P}(\gen12,\gen21)/{\sf I}_{\sf B},
$
where ${\sf I}_{\sf B}$ denotes the ideal generated by
\[
\ZbbZb - \bZbbZ,\
\ZvvZv - \vZvvZ\  \mbox { and }\ 
\dvojiteypsilon - \motylek \hskip .2em.
\]
In the above display we denoted
\begin{eqnarray*}
&\ZbbZb := \gen12(\gen12 \ot \id),\
\bZbbZ := \gen12(\id \ot \gen12),\
\ZvvZv := (\gen21 \ot \id)\gen21,\
\vZvvZ := (\id \ot \gen21)\gen21,&
\\
&\dvojiteypsilon  := \gen21 \circ \gen12
\mbox {\hskip .5em and \hskip .6em}
\motylek := (\gen12 \ot \gen12)\circ \sigma(2,2) \circ  
(\gen21 \ot \gen21),&
\end{eqnarray*}
where $\sigma(2,2)\in \Sigma_4$ is the permutation 
\[
\sigma(2,2) = 
\left(
\begin{array}{cccc}
1 & 2 & 3 & 4
\\
1 & 3 & 2 & 4
\end{array}
\right).
\]

As we argued in~\cite{markl:ba}, the \PROP\ $\sfB$ can be interpreted
as a perturbation of the \PROP\ $\sfhB = F(\sfhb)$ for 
$\frac12$bialgebras mentioned in Example~\ref{eli}. 
More precisely, let $\epsilon$ be a formal parameter,
${\sf I}_{\sf B}^\epsilon$ be the ideal generated by
\[
\ZbbZb - \bZbbZ,\
\ZvvZv - \vZvvZ\  \mbox { and }\ 
\dvojiteypsilon - \epsilon \hskip 1mm \motylek \hskip .2em
\]
and $\sfB_\epsilon := \Gamma_{\tt P}(\gen12,\gen21)/{\sf I}_{\sf B}^\epsilon$.
Then $\sfB_\epsilon$ is a one-dimensional family of deformations of
$\sfhB = \sfB_0$ whose specialization (value) at $\epsilon = 1$ is $\sfB$. 
Therefore, every minimal
model for $\sfB$ can be expected to be a perturbation of the minimal
model for $\sfhB$ described in the following

\begin{theorem}[\cite{markl:ba}]
\label{pozitri_odletam}
The dg free \PROP\
\begin{equation}
\label{41}
({\sf M},\pa_0) = (\Gamma_{\tt P}(\Xi),\pa_0),
\end{equation}
where the generators $\Xi =  \span(\{\xi^n_m\}_{m,n \in I})$ are
as in Example~\ref{Orlik} and the differential $\pa_0$ is given by
formula~(\ref{uz_mi_konci_pobyt_v_Minnesote}), is a 
minimal model of the \PROP\ $\sfhB$ for $\frac12$bialgebras.
\end{theorem}

\noindent 
{\it Proof.}  Clearly, $({\sf M},\pa_0) = 
F(\Omega_{\frac12\tt P}(\Associative * \Associative^\dagger))$. 
The theorem now follows
from Theorem~\ref{eli2} (see also Example~\ref{Orlik}) and from the
fact that the functor $F$ preserves homology isomorphisms, see
Corollary~\ref{8}.%
\qed

The methods developed in this paper were used in~\cite{markl:ba} to prove:

\begin{theorem}
\label{39}
There exists a minimal model $({\sf M},\pa)$ of the \PROP\ \sfB\ for
bialgebras which is a perturbation of the minimal model $({\sf
M},\pa_0)$ of the \PROP\ $\sfhB$ for $\frac12$bialgebras described in
Theorem~\ref{pozitri_odletam},
\[
({\sf M},\pa)  = (\Gamma_{\tt P}(\Xi),\pa_0 + \papert)
\]
for some perturbation $\papert$ which raises the genus and preserves the
path grading. 
\end{theorem}

\noindent 
{\it Proof.}
As shown in~\cite{markl:ba}, a perturbation $\papert$ can be
constructed using standard methods of the homological perturbation
theory because we know, by Theorem~\ref{pozitri_odletam}, 
that $\Gamma_{\tt P}(\Xi)$ is $\pa_0$-acyclic in positive degrees. 
The main problem was to show that
the procedure converges. This was achieved by finding a subspace $X
\subset \Gamma_{\tt P}(\Xi)$ of {\em special elements\/} whose pieces $X(m,n)$ 
satisfy the conditions that:
\begin{itemize}
\item[(i)] \
each $X(m,n)$ is a finite dimensional space spanned by $G$-monomials
with connected~$G$,
\item[(ii)] \
each $X(m,n)$ is $\pa_0$-closed and $\pa_0$-acyclic in positive degrees, 
\item[(iii)] \ each $X(m,n)$ is closed under vertex
insertion (see below) and
\item[(iv)] 
\vskip -2,5mm
\ both $\dvojiteypsilon$ and $\motylek$ belong to $X(2,2)$.
\end{itemize}

Item (iii) means that $X$ is stable under all derivations (not
necessarily differentials) $\omega$ of 
$\Gamma_{\tt P}(\Xi)$ such that $\omega(\Xi) \subset X$. 
The perturbation problem was then solved in $X$ instead of~%
$\Gamma_{\tt P}(\Xi)$. 
It remained to use, in an obvious way,
Corollary~\ref{37}(i) to prove that the object we constructed is
really a minimal model of $\sfB$.
\qed

\vskip 3mm
\noindent
{\bf Dioperads.} In this part we prove that the cobar duals of dioperads
with a replacement rule induce, via functor 
$F_1: \catDiop \to \catPROP$ introduced in Section~\ref{definitions}, 
minimal models in the category of 
\PROP{s}.  Since we are unable to prove the exactness of 
$F_1$, we will need to show first that these models are perturbations
of minimal models of quadratic Koszul \hPROP{s} and then use 
Corollary~\ref{37}(ii).
This approach applies to main examples of~\cite{gan}, i.e.~Lie bialgebras and
infinitesimal bialgebras.

Let ${\calP}$ and $\calQ$ be quadratic operads, with presentations ${\calP} =
\Gamma_{\tt Op}(F)/(R)$ and $\calQ = \Gamma_{\tt Op}(G)/(S)$. We will
consider dioperads created from ${\calP}$ and $\calQ$ by a {\em dioperadic
replacement rule\/}. By this we mean the following. 

As in Example~\ref{posledni_den}, interpret $F$, $G$, $R$ and $S$ as
bicollections. We already observed in Section~\ref{42} that
\[
\Gamma_{\tt D}(F,G)(2,2) \cong
\Gamma_{\frac12\tt P}(F,G)(2,2) \oplus {\rm Ind}^{\Sigma_2 \times
\Sigma_2}_{\{1\}} (F \ot G) \cong G \circ F \oplus   {\rm Ind}^{\Sigma_2 \times
\Sigma_2}_{\{1\}} (F \ot G),
\]
see also~\cite[Section~2.4]{gan} for details. 
The above decomposition is in fact a decomposition of 
$\Gamma_{\tt D}(F,G)(2,2)$ into $\pth$-homogeneous components, namely
\[
G \circ F = 
\span\{f \in \Gamma_{\tt D}(F,G)(2,2);\ \pth(f) = 4\}
\]
and 
\[
 {\rm Ind}^{\Sigma_2 \times
\Sigma_2}_{\{1\}} (F \ot G) = 
\span\{f \in \Gamma_{\tt D}(F,G)(2,2);\ \pth(f) = 3\}.
\]
Given a $(\Sigma_2,\Sigma_2)$-equivariant map
\begin{equation}
\label{repl}
\lambda :  G \circ F \to   {\rm Ind}^{\Sigma_2 \times
\Sigma_2}_{\{1\}} (F \ot G),
\end{equation}
one might consider a subspace
\[
B = B_{\lambda} := \span\{f - \lambda(f);\ f \in G \circ F\} \subset
\Gamma_{\tt D}(F,G)(2,2)
\]
and a quadratic dioperad
\begin{equation}
\label{treti_den_zpatky}
D_\lambda:= \Gamma_{\tt D}(F,G)/(R,B_\lambda,S).
\end{equation}

We say that the map $\lambda$ in~(\ref{repl}) is a {\em replacement
rule}~\cite[Definition~11.3]{fox-markl:ContM97}, 
if it is coherent in the sense that it 
extends to a {\em mixed distributive law\/} between 
operads ${\calP}$ and $\calQ$, 
see~\cite[Section~11]{fox-markl:ContM97} for details. 
An equivalent way to
express this coherence is to say that $D_\lambda$ and 
$F_2({\calP} \diamond \calQ^\dagger)$
are isomorphic as bicollections or, in the terminology
of~\cite[Proposition~5.9]{gan}, that 
$D_\lambda \cong {\calP} \Box \calQ^{\rm op}$, 
see Proposition~\ref{stavka_pilotu}.

\begin{example}
\label{78}
{\rm
An important example is given by an {\em infinitesimal bialgebra\/}
(which we called in~\cite[Example~11.7]{fox-markl:ContM97} {\em a mock
bialgebra\/}). It is a vector space $V$ together with an associative
multiplication $\cdot : V \ot V \to V$ and a coassociative
comultiplication $\Delta : V \to V \ot V$
such that 
\[
\Delta(a \cdot b) = \sum\left(
a_{(1)} \ot a_{(2)}\cdot b + a \cdot b_{(1)} \ot b_{(2)}   
\right)
\]
for any $a,b \in V$.
 
The dioperad $\itIB$ describing infinitesimal bialgebras is given by
$
\itIB = \Gamma_{\tt D}(\gen12,\gen21)/{\it I}_{\itIB},
$
where ${\it I}_{\itIB}$ denotes the dioperadic ideal generated by
\[
\ZbbZb - \bZbbZ,\
\ZvvZv - \vZvvZ\  \mbox { and }\ 
\dvojiteypsilon - \levaplastev  - \pravaplastev\hskip 1mm. 
\]
The dioperad $\itIB$ is created from two copies of the operad 
$\Ass$ for associative algebras using a replacement rule given by
\[
\lambda(\dvojiteypsilon) := \levaplastev  + 
\pravaplastev\hskip 1mm,
\]
see~\cite[Example~11.7]{fox-markl:ContM97} for details.
As before, one may consider a one parameter family 
$\itIB_\epsilon : = \Gamma_{\tt D}(\gen12,\gen21)/{\it I}_{\itIB}^\epsilon$,
where  ${\it I}_{\itIB}^\epsilon$ is the dioperadic ideal generated by
\[
\ZbbZb - \bZbbZ,\
\ZvvZv - \vZvvZ\  \mbox { and }\ 
\dvojiteypsilon - \epsilon \left(\levaplastev  + 
\pravaplastev\hskip 1mm \right)
\]
given by the one parameter
family of replacement rules
\[
\lambda_\epsilon(\dvojiteypsilon) :=  \epsilon \left(\levaplastev  + 
\pravaplastev\hskip 1mm\right).
\]
Let $\sfIB := F_1(\itIB)$ be the \PROP\ generated by the dioperad $\itIB$.
It follows from the above remarks  
that $\sfIB$ is another perturbation of the \PROP\ $\sfhB$
for $\frac12$bialgebras.
}
\end{example}

\begin{example}
\label{80}
{\rm
Recall that a {\em Lie bialgebra\/} is a vector space $V$, with a Lie
algebra structure $[-,-] : V \ot V \to V$ and a Lie diagonal $\delta : V
\to V \ot V$. As in Example~\ref{79}
we assume that the bracket $[-,-]$ is antisymmetric and satisfies
the Jacobi equation and that $\delta$ satisfies the obvious duals of
these conditions, but this time $[-,-]$ and $\delta$ are related by
\[
\delta[a,b] = \sum \left( [a_{(1)},b] \ot a_{(2)} +  [a,b_{(1)}] \ot b_{(2)}
               + a_{(1)} \ot [a_{(2)},b] + b_{(1)} \ot [a,b_{(2)}]  \right)
\]
for any $a,b \in V$, where we used, as usual, the Sweedler notation
$\delta a = \sum a_{(1)} \ot a_{(2)}$ and $\delta b = \sum b_{(1)}
\ot b_{(2)}$.

The dioperad $\itLieB$ for Lie bialgebras is given by
$
\itLieB = \Gamma_{\tt D}(\gen12,\gen21)/{\it I}_{\itLieB},
$
where $\gen 12$ and $\gen 21$ are now {\em antisymmetric\/} generators
and  ${\it I}_{\itLieB}$ denotes the ideal generated by

\vglue 3pt
\[
\Jac 123 + \Jac 231 + \Jac 312 \hskip 1mm , \
\coJac 123 +  \coJac 231 + \coJac 312 \hskip 1mm   \mbox { and }
\hskip 2mm 
\dvojiteypsilonvetsi1212 - 
\levaplastevvetsi1212  - 
\pravaplastevvetsi1212 +
\levaplastevvetsi1221 + \pravaplastevvetsi1221 \hskip 1mm, 
\]
with labels indicating, in the obvious way, the corresponding
permutations of the inputs and outputs.
The dioperad $\itLieB$ is a combination of two  copies of the operad
$\Lie$ for Lie algebras, with the replacement rule

\vskip 2mm 
\[
\lambda\left( \hskip 3pt\dvojiteypsilonvetsi 1212 \hskip 2pt \right) 
:= \levaplastevvetsi1212  +
\pravaplastevvetsi1212 -
\levaplastevvetsi1221 - \pravaplastevvetsi1221 \hskip 1mm, 
\]
see~\cite[Example~11.6]{fox-markl:ContM97}.
One may obtain, as in Example~\ref{78}, a one parameter family
$\itLieB_\epsilon$ of dioperads generated by a one parameter family 
$\lambda_\epsilon$ of 
replacement rules such that $\itLieB_1 = \itLieB$ and $\itLieB_0 =
\ithLieB : = F_2(\sfhlieb)$, where $\sfhlieb$ is the \hPROP\
for $\frac12$Lie bialgebras introduced in Example~\ref{79}. 
Thus, the \PROP\ $\sfLieB :=
F_1 (\itLieB)$ is a perturbation of the \PROP\ $\sfhLieB$ governing
$\frac12$Lie bialgebras. 
}
\end{example}

Examples~\ref{78} and~\ref{80} can be generalized as follows. Each
replacement rule $\lambda$ as in~(\ref{repl}) 
can be extended to a one parameter
family of replacement rules by defining $\lambda_\epsilon :=
\epsilon \cdot \lambda$. This gives a one parameter family
$D_\epsilon : = D_{\lambda_\epsilon}$ of dioperads such that 
$D_1 = D_\lambda$ and $D_0 = {\calP} \diamond \calQ^\dagger$. Therefore
$D_\lambda$ is a perturbation of the dioperad generated by the
\hPROP\ ${\calP} \diamond \calQ^\dagger$. 
This suggests that every minimal model of the \PROP\
$F_1(D_\lambda)$ is a perturbation of a minimal model for 
$F_2({\calP} \diamond \calQ^\dagger)$ which is, as we already know 
from Section~\ref{42}, given by 
$F_2(\Omega_{\frac12\tt P}(({\calP} \diamond \calQ^\dagger)^!)) =
F_2(\Omega_{\frac12\tt P}({\calP}^! *(\calQ^!)^\dagger))$.
The rest of this section makes this idea precise.  

For any quadratic dioperad $D$, there is an obvious candidate for
a minimal model of the \PROP\ $F_1(D)$ generated by $D$, namely the
dg \PROP\  $\Omega_{\tt P}(D^!) = (\Omega_{\tt P}(D^!),\pa) :=
F_1((\Omega_{\tt D}(D^!),\pa))$ generated by the
dioperadic cobar dual $\Omega_{\tt D}(D^!) = 
(\Omega_{\tt D}(D^!),\pa)$ of $D^!$.

The following proposition, roughly speaking, says that the dioperadic
cobar dual of $D_\lambda$ is a perturbation of the cobar dual of the
\hPROP\ $({\calP} \diamond \calQ^\dagger)^! = {\calP}^! *
(\calQ^!)^\dagger$.

\begin{proposition}
\label{44}
Let $D = D_\lambda$ be a dioperad constructed 
from Koszul quadratic operads ${\calP}$
and $\calQ$ using a replacement rule $\lambda$.
Consider the canonical decomposition 
\[
(\Omega_{\tt D}(D^!),\pa_0 + \pa_\pth)
\]
of the differential in the dioperadic bar
construction $(\Omega_{\tt D}(D^!),\pa)$. Then 
\begin{equation}
\label{801}
(\Omega_{\tt D}(D^!),\pa_0) \cong F_2(\Omega_{\frac12\tt P}({\calP}^! * 
(\calQ^!)^\dagger)).
\end{equation}
\end{proposition}

\noindent 
{\it Proof.} We already observed that, in the terminology
of~\cite{gan}, $D \cong {\calP} \Box \calQ^{\rm op}$. This implies, 
by~\cite[Proposition~5.9(b)]{gan}, 
that $D^! \cong (\calQ^!)^{\rm op} \Box {\calP}^!$
which clearly coincides, as a bicollection, with our 
${\calP}^! * (\calQ^!)^\dagger$. The rest of the proposition
follows from the description of $D^!$ given in~\cite{gan},
the behavior of the replacement rule $\lambda$ with respect to the path
grading, and definitions.%
\qed

\begin{remark}
{\rm\ 
Since, as a non-differential dioperad, 
$\Omega_{\tt D}(D) = \Lambda^{-1} \Gamma_{\tt D} 
(\uparrow {\bar D}^*)$, where $\uparrow$ denotes 
the suspension of a graded bicollection, 
$\Lambda^{-1}$ the sheared desuspension of a
dioperad and $\bar D^*$ the linear dual of the 
augmentation ideal of $D$, see Sections 1.4, 2.3, and 3.1 
of~\cite{gan} for details, the \PROP\ $(\Omega_{\tt P} (D),\pa)$
may be constructed from scratch as $\Omega_{\tt P}(D^!) = \Lambda^{-1}
\Gamma_{\tt P} (\uparrow \bar D^*)$ 
with a differential coming from the ``vertex expansion''
(also called edge insertion). Thus, the \PROP\ $(\Omega_{\tt P}
(D),\pa)$ may be thought of as a \emph{naive cobar dual} of
$F_1(D)$, as opposed to the categorical cobar dual~\cite[Section
4.1.14]{ginzburg-kapranov:DMJ94}.

Perhaps, one can successfully develop quadratic and Koszul duality
theory for \PROP{s}, using this naive cobar dual by analogy 
with~\cite{ginzburg-kapranov:DMJ94,gan}. We are reluctant to emphasize
 $(\Omega_{\tt P} (D),\pa)$ as a PROP cobar dual of the \PROP\
$F_1(D)$, because we do not know how this naive cobar dual is 
related to the categorical one.
}
\end{remark}

The following theorem generalizes a result of
Kontsevich~\cite{kontsevich:message} for $D = \itLieB$.

\begin{theorem}
\label{45}
Under the assumptions of Proposition~\ref{44},
$(\Omega_{\tt P}(D^!),\pa)$ is a minimal model of the \PROP\ $F_1(D)$.
\end{theorem}

\noindent 
{\it Proof of Theorem~\ref{45}.}
We are going to use Corollary~\ref{37}(ii).
It is straightforward to verify that $H_0((\Omega_{\tt
{\calP}}(D^!),\pa) \cong F_1(D)$. Equation~(\ref{801}) gives
\[
\Omega_{\tt P}(D^!) \cong F(\Omega_{\frac12\tt P}({\calP}^! * 
(\calQ^!)^\dagger)), 
\]
therefore the $\pa_0$-acyclicity 
of $\Omega_{\tt P}(D^!)$ follows from the exactness of the
functor $F$ stated in Theorem~\ref{polynomial}.%
\qed

\begin{example}
\label{390}
{\rm\
By Theorem~\ref{45}, the dg \PROP\  $\Omega_{\tt P}({\it IB}^!)$, where
the quadratic dual ${\itIB}^!$ of the dioperad $\itIB$ for
infinitesimal bialgebras is described in~\cite{gan}
as ${\itIB}^! = \Ass^{\rm op} \Box \Ass$, is a 
minimal model of the \PROP\ $\sfIB  = F_1({\it IB})$ 
for infinitesimal bialgebras. The dg \PROP\
$\Omega_{\tt P}({\it IB}^!)$ has a
form $(\Gamma_{\tt P}(\Xi),\pa_0 + \pa_\pth)$, where $\Xi$ and $\pa_0$
are the same as in Example~\ref{Orlik}. The path part $\pa_\pth$ of the
differential is trivial on generators $\xi^m_n$ with $m+n \leq 4$,
therefore the easiest example of the path part is provided by

\[
\pa(\dvatri) = \pa_0(\dvatri) + \jitka + \jitkainv + \anna - \annainv
\]
where
\[
\pa_0(\dvatri) = 
\dvacarkatri - \dvaZbbZb + \dvabZbbZ
\]
is the same as in Example~\ref{Orlik}. We encourage the reader to
verify that
\begin{eqnarray*}
&
\pth(\dvatri) = \pth(\dvacarkatri) = \pth(\dvaZbbZb) = \pth(\dvabZbbZ)
= 6,
&
\\
&
\pth(\jitka) = \pth(\jitkainv) = 5\
\mbox { and }\
\pth(\anna) = \pth(\annainv) = 4.
&
\end{eqnarray*}

Similarly, the dg \PROP\ $\Omega_{\tt P}({\it LieB}^!)
= \Omega_{\tt P}(\Com^{\rm op} \Box \Com)$ is a minimal model of the 
\PROP\ $\sfLieB := F_1({\it LieB})$ for
Lie bialgebras.
}
\end{example}

\section{Classical graph cohomology}
\label{classical}

Here we will reinterpret minimal models for the Lie bialgebra
\PROP\ $\sfLieB = F_1 (\itLieB)$ and the infinitesimal bialgebra
\PROP\ $\sfIB = F_1 (\itIB)$ given by Theorem~\ref{45} and
Example~\ref{390} as graph complexes.  \smallskip

\noindent
{\bf The commutative case}. Consider the set of connected
$(m,n)$-graphs $G$ for $m,n \ge 1$ in the sense of
Section~\ref{fprops}. An \emph{orientation} on an $(m,n)$-graph $G$ is
an orientation on $\nr^{v(G)} \oplus \nr^{m} \oplus \nr^n$, i.e., the
choice of an element in $\det \nr^{v(G)} \otimes \det \nr^m
\otimes \det \nr^n$ up to multiplication by a positive real number.
This is equivalent to an orientation on $\nr^{e(G)} \oplus H_1
(\abs{G}; \nr)$, where $e(G)$ is the set of (all) edges of $G$; to
verify this, consider the cellular chain complex of the geometric
realization $\abs{G}$ , see for example \cite[Proposition
B.1]{thurston:math.QA/9901110} and \cite[Proposition
5.65]{markl-shnider-stasheff:book}.

\begin{figure}[t]
\begin{center}
\unitlength=.6cm
\begin{picture}(5,7)(0.00,-3)
\thicklines
\put(0.00,0){\makebox(0.00,0.00){$\bullet$}}
\put(2.00,-2){\makebox(0.00,0.00){$\bullet$}}
\put(2.00,2){\makebox(0.00,0.00){$\bullet$}}
\put(2.00,3){\makebox(0.00,0.00){$\bullet$}}
\put(4.00,0){\makebox(0.00,0.00){$\bullet$}}
\put(2.00,-3.00){\line(0,1){1}}
\put(2.00,2.00){\line(0,1){2}}
\put(2.00,-2.00){\line(1,1){2}}
\put(2.00,-2.00){\line(-1,1){2}}
\put(0.00,0){\line(1,1){2}}
\put(4.00,0){\line(-1,1){2}}
\put(0.00,0){\line(2,3){2}}
\put(4.00,0){\line(-2,3){2}}
\end{picture}
\end{center}
\caption{\label{aut-zero}
A graph vanishing in the quotient by the automorphism group.}
\end{figure}
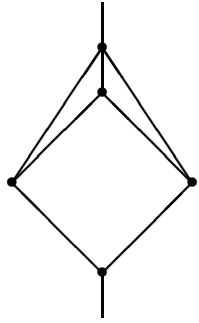

Thus, an orientation on a connected $(m,n)$-graph $G$ is equivalently
given by an ordering of the set $e(G)$ along with the choice of an
orientation on $H_1 (\abs{G}; \nr)$ up to permutations and changes of
orientation on $H_1 (\abs{G}; \nr)$ of even total parity. Consider the
set of isomorphism classes of oriented $(m,n)$-graphs and take its
$k$-linear span. More precisely, we should rather speak about a
colimit with respect to graph isomorphisms, as in
Section~\ref{fprops}. In particular, if a graph $G$ admits an
orientation-reversing automorphism, such as the graph in
Figure~\ref{aut-zero}, then $G$ gets identified with $G^-$, which will
vanish after passing to the following quotient.
Let $G(m,n)$ be the
quotient of this space by the subspace spanned by
\[
G + G^- \qquad \mbox{ for each oriented graph } G,
\]
where $G^-$ is the same graph as $G$, taken with the opposite
orientation. Each space $G(m,n)$ is bigraded by the genus and the
number of \emph{interior edges} (i.e., edges other than legs) of the
graph. Let $G_g^q = G_g^q(m,n)$ denote the subspace spanned by graphs
of genus $g$ with $q$ interior edges for $g, q \ge 0$. Computing the
Euler characteristic of $\abs{G}$ in two ways, we get an identity
$\abs{v(G)} - q = 1 - g$. A graph $G \in G_g (m,n)$ has a maximal
number of interior edges, if each vertex of $G$ is trivalent, in which
case we have $3 \abs{v(G)} = 2q + m + n$, whence $q = 3g - 3 + m + n$
is the top degree in which $G_g^q(m,n) \ne 0$.

Define a differential
\[
\pa: G_g^q \to G_g^{q+1},
\]
so that $\pa^2 =0$, as follows:
\[
\pa G := \sum_{\{G' \; | \; G'/e = G\}} G',
\]
where the sum is over the isomorphism classes of connected
$(m,n)$-graphs $G'$ whose contraction along an edge $e \in e(G')$ is
isomorphic to $G$. We will induce an orientation on $G'$ by first
choosing an ordering of the set of edges of $G$ and an orientation on
$H_1 (\abs{G}; \nr)$ in a way compatible with the orientation of
$G$. Then we will append the edge $e$ which is being contracted at the
end of the list of the edges of $G$. Since we have a canonical
isomorphism $H_1 (\abs{G'}; \nr) \stackrel{\sim}{\to} H_1 (\abs{G};
\nr)$, an orientation on the last space induces one on the first. This
gives an orientation on $G'$. An example is given below.
\[
\pa \left(\hskip 3pt
{\unitlength=.46pt
\begin{picture}(24.00,30.00)(0.00,3.00)
\bezier{20}(10.00,10.00)(15,5)(20.00,0.00)
\bezier{20}(10.00,10.00)(5,5)(0.00,0.00)
\bezier{20}(10.00,10.00)(15,15)(20.00,20.00)
\bezier{20}(0.00,20.00)(5,15)(10.00,10.00)
\put(0,-5){\makebox(0,0)[t]{\scriptsize 1}}
\put(20,-5){\makebox(0,0)[t]{\scriptsize 2}}
\put(0,25){\makebox(0,0)[b]{\scriptsize 1}}
\put(20,25){\makebox(0,0)[b]{\scriptsize 2}}
\end{picture}} \hskip 2pt
\right) = 
\dvojiteypsilonvetsi1212 - 
\levaplastevvetsi1212  - 
\pravaplastevvetsi1212 +
\levaplastevvetsi1221 +
\pravaplastevvetsi1221
\]
In this figure we have oriented graphs, which are provided with a
certain canonical orientation that may be read off from the picture.
The rule of thumb is as follows. An \emph{orientation on the
composition} of two graphs is given by (1) reordering the edges of the
first, lower, graph in such a way that the output legs follow the
remaining edges, (2) reordering the edges of the second, upper, graph
in such a way that the input legs precede the remaining edges, and (3)
after grafting, putting the edges of the second graph after the edges
of the first graph. The resulting ordering should look like this: the
newly grafted edges in the middle, preceded by the remaining edges of
the first graph and followed by the remaining edges of the second
graph.  We remind the reader that we place the inputs at the bottom of
a graph and the outputs on the top.

\begin{theorem}
\label{comm-graph}
The graph complex in the commutative case is acyclic everywhere but at
the top term $G_g^{3g-3 +m +n}$. The graph cohomology can be computed
as follows.
\[
H^q (G_g^* (m,n), \pa) =
\left\{
\begin{array}{ll}
{\sfLieB}_g^0 (m,n) & \mbox{for $q = 3g -3 + m + n$},\\
0 & \mbox{otherwise},
\end{array}
\right.
\]
where $\sfLieB_g^0 (m,n)$ is the subspace of the $(m,n)$th component
of the Lie bialgebra \PROP\ $\sfLieB = F_1 (\itLieB)$ consisting of
linear combinations of connected graphs of genus $g$, see the
presentation of the corresponding dioperad $\itLieB$ in
Example~\ref{80}.
\end{theorem}

\begin{remark}
{\rm\ 
The acyclicity of the graph complex $G^*_g (m,n)$ has been proven in
Kontsevich's message \cite{kontsevich:message}, whose method we have
essentially used in this paper.
}
\end{remark}

\noindent 
{\it Proof}. The dioperad $\itLieB$ may be represented as a $\Box$
product of the Lie operad $\Lie$ and the Lie co-operad $\Lie^{\rm
op}$: $\itLieB = \Lie \Box \Lie^{\rm op}$ --- see \cite[Section
5.2]{gan}. The dioperadic quadratic dual $\itLieB^!$ is then
$\Com^{\rm op} \Box \Com$, so that $\itLieB^!(m,n) \cong k$ with a
trivial action of $(\Sigma_n,\Sigma_n)$ for each pair $(m,n)$, $m,n
\ge 1$. Then the subcomplex $(\Omega^0_{\tt P}({\itLieB}^!),\pa)
\subset (\Omega_{\tt P}({\itLieB}^!),\pa)$ spanned by connected graphs
is isomorphic to the graph complex $(G^* (m,n), \pa)$. Now the result
follows from Theorem~\ref{45}.  \qed
\smallskip

\noindent
{\bf The associative case}. Consider connected, oriented
$(m,n)$-graphs $G$ for $m,n \ge 1$, as above, now with a \emph{ribbon}
structure at each vertex, by which we mean orderings of the set
$\In(v)$ of incoming edges and the set $\Out(v)$ of outgoing edges at
each vertex $v \in v(G)$. It is convenient to think of an equivalent
cyclic ordering (i.e., ordering up to cyclic permutation) of the set
$e(v) = \In(v) \cup \Out(v)$ of all the edges incident to a vertex $v$
in a way that elements of $\In(v)$ precede those of $\Out(v)$. Let
$\RG (m,n)$ be the linear span of isomorphism classes of connected
oriented ribbon $(m,n)$-graphs modulo the relation $G + G^- = 0$, with
$\RG^q_g (m,n)$ denoting the subspace of graphs of genus $g$ with $q$
interior edges.  The same formula
\[
\pa G := \sum_{\{G' \; | \; G'/e = G\}} G'
\]
defines a differential, except that in the ribbon case, when we
contract an edge $e \in e(G')$, we induce a cyclic ordering on the set
of edges adjacent to the resulting vertex by an obvious operation of
insertion of the ordered set of edges adjacent to the edge $e$ through
one of its vertices into the ordered set of edges adjacent to $e$
through its other vertex. An orientation is induced on $G'$ in the
same way as in the commutative case. An example is shown in the
following display.
\[
\pa(\dvactyri) =
{
\unitlength=.1pt
\begin{picture}(96.00,120.00)(-8.00,0.00)
\qbezier(40.00,60.00)(20.00,80.00)(0.00,100.00)
\qbezier(40.00,60.00)(60.00,80.00)(80.00,100.00)
\put(40.00,60.00){\line(0,-1){20.00}}
\qbezier(40.50,39.00)(49.00,14.50)(53.50,0.00)
\qbezier(40.00,40.00)(32.50,16.50)(26.00,0.00)
\qbezier(40.00,40.00)(60.00,20.00)(80.00,0.00)
\qbezier(40.00,40.00)(20.00,20.00)(0.00,0.00)
\end{picture}}
-
{\unitlength=.3pt
\begin{picture}(40,20)(-34,0)
\put(0,0){\line(0,1){10}}
\put(-10,20){\line(0,1){10}}
\put(-30,20){\line(0,1){10}}
\put(-10,0){\bezier{40}(0.00,0.00)(-10.00,10.00)(-20.00,20.00)}
\put(-10,0){\bezier{40}(0.00,20.00)(-10.00,10.00)(-20.00,0.00)}
\put(0,10){\bezier{20}(0.00,0.00)(-5.00,5.00)(-10.00,10.00)}
\put(-19,0){\line(0,1){10}}
\end{picture}}
-
{\unitlength=.32pt
\begin{picture}(40,20)(-2,0)
\put(10,0){\line(0,1){30}}
\put(20,0){\line(0,1){10}}
\put(30,20){\line(0,1){10}}
\put(-10,0){\bezier{40}(0.00,0.00)(10.00,10.00)(20.00,20.00)}
\put(0,0){\bezier{20}(0.00,0.00)(5.00,10.00)(10.00,20.00)}
\put(20,10){\bezier{20}(0.00,0.00)(5.00,5.00)(10.00,10.00)}
\put(20,10){\bezier{20}(0.00,0.00)(-5.00,5.00)(-10.00,10.00)}
\end{picture}}
+ \dvabbZbbZ - \dvabZbbZb + \dvaZbbZbb - 
{\unitlength=.3pt
\begin{picture}(40,20)(-2,0)
\put(0,0){\line(0,1){10}}
\put(10,20){\line(0,1){10}}
\put(30,20){\line(0,1){10}}
\put(10,0){\bezier{40}(0.00,0.00)(10.00,10.00)(20.00,20.00)}
\put(10,0){\bezier{40}(0.00,20.00)(10.00,10.00)(20.00,0.00)}
\put(0,10){\bezier{20}(0.00,0.00)(5.00,5.00)(10.00,10.00)}
\put(20,0){\line(0,1){10}}
\end{picture}
}
-
{\unitlength=.32pt
\begin{picture}(40,20)(-34,0)
\put(-10,0){\line(0,1){30}}
\put(-20,0){\line(0,1){10}}
\put(-30,20){\line(0,1){10}}
\put(10,0){\bezier{40}(0.00,0.00)(-10.00,10.00)(-20.00,20.00)}
\put(-20,10){\bezier{20}(0.00,0.00)(-5.00,5.00)(-10.00,10.00)}
\put(-20,10){\bezier{20}(0.00,0.00)(5.00,5.00)(10.00,10.00)}
\put(0,0){\bezier{20}(0.00,0.00)(-5.00,10.00)(-10.00,20.00)}
\end{picture}
}
-
{\unitlength=.32pt
\begin{picture}(40,20)(-34,0)
\put(0,0){\line(0,1){10}}
\put(-10,0){\line(0,1){30}}
\put(-20,0){\line(0,1){10}}
\put(-30,20){\line(0,1){10}}
\put(0,10){\bezier{20}(0.00,0.00)(-5.00,5.00)(-10.00,10.00)}
\put(-20,10){\bezier{20}(0.00,0.00)(-5.00,5.00)(-10.00,10.00)}
\put(-30,0){\bezier{40}(0.00,0.00)(10.00,10.00)(20.00,20.00)}
\end{picture}
}
-
{\unitlength=.32pt
\begin{picture}(40,20)(-2,0)
\put(0,0){\line(0,1){10}}
\put(10,0){\line(0,1){30}}
\put(20,0){\line(0,1){10}}
\put(30,20){\line(0,1){10}}
\put(0,10){\bezier{20}(0.00,0.00)(5.00,5.00)(10.00,10.00)}
\put(20,10){\bezier{20}(0.00,0.00)(5.00,5.00)(10.00,10.00)}
\put(30,0){\bezier{40}(0.00,0.00)(-10.00,10.00)(-20.00,20.00)}
\end{picture}
}
- \dvabZbbbZ - \dvaZbbbZb
\]

A vanishing theorem, see below, also holds in the ribbon-graph case.
The proof is similar to the commutative case: it uses Theorem~\ref{45}
and the fact that $\itIB = \Ass \Box \Ass^{\rm op}$ and $\itIB^! =
\Ass^{\rm op} \Box \Ass$, see Example~\ref{390}.

\begin{theorem}
\label{ribbon-graph}
\begin{sloppypar}
  The ribbon graph complex is acyclic everywhere but at the top term
  $\RG_g^{3g-3 +m +n}$. The ribbon graph cohomology can be computed as
  follows.
\[
H^q (\RG_g^* (m,n), \pa) =
\left\{
\begin{array}{ll}
{\sfIB}_g^0 (m,n) & \mbox{for $q = 3g -3 + m + n$},\\
0 & \mbox{otherwise},
\end{array}
\right.
\]
where $\sfIB_g^0 (m,n)$ is the subspace of the $(m,n)$th component of
the infinitesimal bialgebra \PROP\ $\sfIB = F_1 (\itIB)$ consisting of
linear combinations of connected ribbon graphs of genus $g$, see the
presentation of the corresponding dioperad $\itIB$ in
Example~\ref{78}.
\end{sloppypar}
\end{theorem}

\begin{remark}
  {\rm Note that our notion of the genus is not the same as the one
    coming from the genus of an oriented surface associated to the
    graph, usually used for ribbon graphs. Our genus is just the first
    Betti number of the surface. }
\end{remark}

\bibliographystyle{smfplain}

\begin{thebibliography}{10}

\bibitem{aguiar}
{\scshape M.~Aguiar} -- {\og Infinitesimal {Hopf} algebras\fg}, \emph{New
  trends in Hopf algebra theory (La Falda, 1999)}, Contemporary Math., vol.
  267, Amer. Math. Soc., 2000, p.~1--29.

\bibitem{culler-vogtmann:IM86}
{\scshape M.~Culler {\normalfont \smfandname} K.~Vogtmann} -- {\og Moduli of
  graphs and automorphisms of free groups\fg}, \emph{Invent. Math.} \textbf{84}
  (1986), p.~91--119.

\bibitem{fox-markl:ContM97}
{\scshape T.~Fox {\normalfont \smfandname} M.~Markl} -- {\og Distributive laws,
  bialgebras, and cohomology\fg}, \emph{Operads: Proceedings of Renaissance
  Conferences} (J.-L. Loday, J.~Stasheff {\normalfont \smfandname} A.~Voronov,
  \smfedsname), Contemporary Math., vol. 202, Amer. Math. Soc., 1997,
  p.~167--205.

\bibitem{gan}
{\scshape W.~Gan} -- {\og Koszul duality for dioperads\fg}, \emph{Math. Res.
  Lett.} \textbf{10} (2003), no.~1, p.~109--124.

\bibitem{getzler-kapranov:CompM98}
{\scshape E.~Getzler {\normalfont \smfandname} M.~Kapranov} -- {\og Modular
  operads\fg}, \emph{Compositio Math.} \textbf{110} (1998), no.~1, p.~65--126.

\bibitem{ginzburg-kapranov:DMJ94}
{\scshape V.~Ginzburg {\normalfont \smfandname} M.~Kapranov} -- {\og Koszul
  duality for operads\fg}, \emph{Duke Math. J.} \textbf{76} (1994), no.~1,
  p.~203--272.

\bibitem{kontsevich:93}
{\scshape M.~Kontsevich} -- {\og Formal (non)commutative symplectic
  geometry\fg}, The Gel'fand mathematics seminars 1990--1992, Birkh\"auser,
  1993.

\bibitem{kont:feyn}
\bysame , {\og Feynman diagrams and low-dimensional topology\fg}, \emph{First
  European Congress of Mathematics~{II}}, Progr. Math., vol. 120,
  Birkh\"{a}user, Basel, 1994.

\bibitem{kontsevich:message}
\bysame , {\og An e-mail message to {M.~Markl}\fg}, November 2002.

\bibitem{maclane:RiceUniv.Studies63}
{\scshape S.~{Mac~Lane}} -- {\og Natural associativity and commutativity\fg},
  \emph{Rice Univ. Stud.} \textbf{49} (1963), no.~1, p.~28--46.

\bibitem{madsen-weiss}
{\scshape I.~Madsen {\normalfont \smfandname} M.~S. Weiss} -- {\og The stable
  moduli space of {Riemann} surfaces: {Mumford's} conjecture\fg}, Preprint {\tt
  math.AT/0212321}, December 2002.

\bibitem{markl:JPAA96}
{\scshape M.~Markl} -- {\og Cotangent cohomology of a category and
  deformations\fg}, \emph{J. Pure Appl. Algebra} \textbf{113} (1996), no.~2,
  p.~195--218.

\bibitem{markl:zebrulka}
\bysame , {\og Models for operads\fg}, \emph{Comm. Algebra} \textbf{24} (1996),
  no.~4, p.~1471--1500.

\bibitem{markl:ba}
\bysame , {\og A resolution (minimal model) of the {PROP} for bialgebras\fg},
  \emph{J. Pure Appl. Algebra} \textbf{205} (2006), no.~2, p.~341--374.

\bibitem{markl-shnider-stasheff:book}
{\scshape M.~Markl, S.~Shnider {\normalfont \smfandname} J.~D. Stasheff} --
  \emph{Operads in algebra, topology and physics}, Mathematical Surveys and
  Monographs, vol.~96, American Mathematical Society, Providence, Rhode Island,
  2002.

\bibitem{merkulov:defquant}
{\scshape S.~Merkulov} -- {\og {PROP} profile of deformation quantization\fg},
  Preprint {\tt math.QA/0412257}, December 2004.

\bibitem{penner:JDG88}
{\scshape R.~Penner} -- {\og Perturbative series and the moduli space of
  {Riemann} surfaces\fg}, \emph{J. Differential Geom.} \textbf{27} (1988),
  no.~1, p.~35--53.

\bibitem{thurston:math.QA/9901110}
{\scshape D.~Thurston} -- {\og Integral expressions for the {V}assiliev knot
  invariants\fg}, Preprint {\tt math.QA/9901110}, January 1999.

\bibitem{vallette:thesis}
{\scshape B.~Vallette} -- {\og Dualit{\' e} de {Koszul} des {PROP}s\fg},
  \smfphdthesisname, Universit{\'e} Louis Pasteur, 2003.

\end{thebibliography}
\def\cprime{$'$}
\providecommand{\bysame}{\leavevmode ---\ }
\providecommand{\og}{``}
\providecommand{\fg}{''}
\providecommand{\smfandname}{et}
\providecommand{\smfedsname}{\'eds.}
\providecommand{\smfedname}{\'ed.}
\providecommand{\smfmastersthesisname}{M\'emoire}
\providecommand{\smfphdthesisname}{Th\`ese}

\end{document}